\documentclass[a4paper,USenglish,cleveref, autoref, thm-restate]{lipics-v2021}

\pdfoutput=1 

\usepackage[numbers,sort&compress]{natbib}
\title{Peel-and-Bound: Generating Stronger Relaxed Bounds with Multivalued Decision Diagrams}
\titlerunning{Peel-and-Bound: Generating Stronger Bounds with Decision Diagrams}

\author{Isaac Rudich}{Mathematics and Industrial Engineering Department, Polytechnique Montréal, Canada }{isaac.rudich@polymtl.ca}{}{}

\author{Quentin Cappart}{Computer Engineering and Software Engineering Department, Polytechnique Montréal, Canada}{quentin.cappart@polymtl.ca}{}{}

\author{Louis-Martin Rousseau}{Mathematics and Industrial Engineering Department, Polytechnique Montréal, Canada \and \url{https://hanalog.ca/person/louis-martin-rousseau/} }{louis-martin.rousseau@polymtl.ca}{}{}

\authorrunning{I. Rudich and Q. Cappart and L.M. Rousseau}
\Copyright{Isaac Rudich and Louis-Martin Rousseau and Quentin Cappart}

\category{} 
\relatedversion{} 

\supplement{}
\supplementdetails[linktext={https://github.com/IsaacRudich/PnB$\_$SOP}]{Source Code for Experiments }{https://github.com/IsaacRudich/PnB_SOP} 

\nolinenumbers 

\EventEditors{Christine Solnon}
\EventNoEds{1}
\EventLongTitle{28th International Conference on Principles and Practice of Constraint Programming (CP 2022)}
\EventShortTitle{CP 2022}
\EventAcronym{CP}
\EventYear{2022}
\EventDate{July 31--August 8, 2022}
\EventLocation{Haifa, Israel}
\EventLogo{}
\SeriesVolume{235}
\ArticleNo{15}
\usepackage[ruled,vlined, linesnumbered]{algorithm2e}

\begin{document}
\definecolor{safeA2}{RGB}{212,17,89}

\maketitle

\begin{abstract}
Decision diagrams are an increasingly important tool in cutting-edge solvers for discrete optimization. However, the field of decision diagrams is relatively new, and is still incorporating the library of techniques that conventional solvers have had decades to build. We drew inspiration from the \emph{warm-start} technique used in conventional solvers to address one of the major challenges faced by decision diagram based methods. Decision diagrams become more useful the wider they are allowed to be, but also become more costly to generate, especially with large numbers of variables. We present a method of \emph{peeling} off a sub-graph of previously constructed diagrams and using it as the initial diagram for subsequent iterations that we call \emph{peel-and-bound}. We test the method on the \emph{sequence ordering problem}, and our results indicate that our \emph{peel-and-bound} scheme generates stronger bounds than a branch-and-bound scheme using the same propagators, and at significantly less computational cost.
\end{abstract}

\section{Introduction}
\label{sec:introduction}
Multivalued decision diagrams (MDDs) are a useful graphical tool for compactly storing the solution space of discrete optimization problems. In the last few years, a staggering number of new applications for MDDs have been proposed~\cite{ddsurvey}, such as representing global constraints~\cite{global1,global2,global3}, handling stochastic variables~\cite{stochastic1,stochastic2}, and performing post-optimality analysis~\cite{postopt}. MDDs are particularly useful for generating strong dual bounds~\cite{cappart2019improving,dual2,dual1,dual3}, especially on optimization problems where linear relaxations perform poorly. There is a subset of MDD research that uses a highly paralellizable \emph{branch-and-bound} algorithm based on decision diagrams~\cite{BnB,ddo,bnb3,parjadis2021improving} to maximize the utility of using MDD based relaxations. This paper furthers the work on the decision diagram based branch-and-bound by introducing a method, referred to as \textit{peel-and-bound}, of reusing the graphs generated at each iteration of the algorithm. Specifically, the contributions are as follows: (1) we present the \emph{peel-and-bound} algorithm, (2) we identify several heuristic decisions that can be used to adjust \emph{peel-and-bound}, and discuss their implications, (3) we show that \emph{peel-and-bound} outperforms \emph{branch-and-bound} on the \emph{sequence ordering problem} (SOP), and (4) we provide insight into how the algorithm can be applied to other problems.

The paper is structured as follows. The next section provides the necessary technical background information and notation, as well as implementation details for the decision diagram relaxations used in our experiments. In Section \ref{sec:pnb} we introduce the core contribution, namely the \emph{peel-and-bound} procedure. The algorithm is presented, and its limitations are discussed. Computational experiments are proposed and discussed in Section \ref{sec:results}. 

\section{Technical Background}
\label{sec:background}

The idea of using multivalued decision diagrams (MDDs) to generate relaxed bounds for optimization problems was introduced by Andersen et al.~(2010)~\cite{andersen2007}. This has been generalized by Hadzic et al.~(2008)~\cite{hadzic2008} and Hoda et al.~(2010)~\cite{hoda2010}. 
Following those papers, Bergman et al.~\cite{BnB,BnB2} demonstrated the potential for a decision diagram based branch-and-bound solver to be effective, and provided an efficient parallelization scheme. Gillard et al.~(2021)~\cite{ddo} further improved the decision diagram based branch-and-bound solver by adding pruning techniques that can be used while the decision diagrams are being constructed, as well as to remove nodes from the branch-and-bound queue. 

This paper presents a new peel-and-bound scheme for combining restricted and relaxed decision diagrams to find exact solutions. This section provides the required technical background on  how decision diagrams can be used to model sequencing problems, and how to construct restricted/relaxed diagrams. It also introduces the notations used in this paper, and details the existing algorithms considered in our experiments. 
 
\subsection{Decision Diagrams (DDs)}
\label{sec:background:dd}
Let $\mathcal{P}$ be an instance of a discrete minimization problem with $n$ variables $\{x_1,...,x_n\}$, let $Sol(\mathcal{P})$ be the set of feasible solutions to $\mathcal{P}$, let $z^*(\mathcal{P})$ be an optimal solution to $\mathcal{P}$, and let $D(x_i)$ be the domain of variable $x_i, i \in \{1, ..., n\}$. Let $\mathcal{M}$ be a multivalued decision diagram that contains potential solutions to $\mathcal{P}$. $\mathcal{M}$ is a directed acyclic graph divided into $n+1$ layers; let $\ell_u$ be the index of the layer containing node $u, u \in \mathcal{M}$, and let $L_i$ be the set containing the nodes on layer $i$. Layer $1$ contains only a root node $r$ (with no \emph{in} arcs), and layer $n+1$ contains just a terminal node $t$ (with no \emph{out} arcs). Each arc $a_{uv} \in \mathcal{M}$ goes from a node $u$ on layer $\ell_u \in \{1,...,n\}$ to a node $v$ on layer $\ell_{u+1}$ ($\ell_{u+1} = \ell_v$). Each arc $a_{uv}$ has a label representing the assignment of variable $x_{\ell_u}$ to $l \in D(x_{\ell_u})$. An arc $a_{uv}$ with label $l$ ($a_{uv} \rightarrow l$) also has a value $v(a_{uv})$ equal to the value of being at node $u$ and assigning $x_{\ell_u}$ to $l$ ($x_{\ell_u} = l$). For simplicity, we sometime refer to $v(a_{uv})$ as $v(a)$. Thus, each path from $r$ to $t$ represents the assignment of the $n$ variables to values, and a potential solution to $\mathcal{P}$.

Let $Sol(\mathcal{M})$ be the set of all paths in $\mathcal{M}$ from $r$ to $t$, and let $T^*(u)$ be the value of the shortest path from $r$ to a node $u$. If $Sol(\mathcal{M}) = Sol(\mathcal{P})$, then $\mathcal{M}$ perfectly represents the solution space of $\mathcal{P}$, and we call $\mathcal{M}$ \emph{exact}. If $\mathcal{M}$ is exact, then the value of the shortest path through the diagram is $z^*(\mathcal{P})$ (an optimal solution to $\mathcal{P}$). Let the shortest path through $\mathcal{M}$ be $z^*(\mathcal{M})$. If $Sol(\mathcal{M}) \subseteq Sol(\mathcal{P})$, then $\mathcal{M}$ represents only feasible solutions to $\mathcal{P}$, but does not necessarily represent all feasible solutions to $\mathcal{P}$. In this case, we call $\mathcal{M}$ \emph{restricted}, and use the notation $\overline{\mathcal{M}}$ to mean that $\mathcal{M}$ is restricted. The shortest path through $\overline{\mathcal{M}}$ is not guaranteed to be optimal, but it is guaranteed to be feasible. If $Sol(\mathcal{P}) \subseteq Sol(\mathcal{M})$, then $\mathcal{M}$ represents all of the feasible solutions to $\mathcal{P}$, but potentially represents infeasible solutions as well. In this case, we call $\mathcal{M}$ \emph{relaxed}, and use the notation $\underline{\mathcal{M}}$ to mean that $\mathcal{M}$ is relaxed. The shortest path through $\underline{\mathcal{M}}$ is guaranteed to be at least as good as $z^*(\mathcal{P})$, but is not guaranteed to be feasible.

 Constructing an exact decision diagram for $\mathcal{P}$ is often intractable for large values of $n$. Observe that having an exact decision diagram means that the solution to $\mathcal{P}$ can be read in polynomial time by recursively calculating the shortest path through $\mathcal{M}$, so creating an exact decision diagram for NP-hard problems, such as for the \emph{travelling salesperson problem} (TSP), is NP-hard as well~\cite{DDforO}. The focus of most research that uses decision diagrams for optimization is on the construction of $\overline{\mathcal{M}}$ and/or $\underline{\mathcal{M}}$. Let $w = w(\mathcal{M})$ be the width of the largest layer of $\mathcal{M}$. The creation of an exact decision diagram potentially leads to $w$ being an exponential function of $n$, but when creating $\overline{\mathcal{M}}$ and/or $\underline{\mathcal{M}}$, $w$ can be constrained to be any natural number, limiting the number of operations construction will take. Let $w_m$ be the largest width allowed during construction. As $w_m$ approaches the width necessary to create an exact decision diagram, $z^*(\overline{\mathcal{M}})$ and $z^*(\underline{\mathcal{M}})$ approach $z^*(\mathcal{P})$, but the number of operations necessary to construct the diagram also increases.
 
\begin{figure}[h!]
    \centering
\includegraphics[scale=.6]{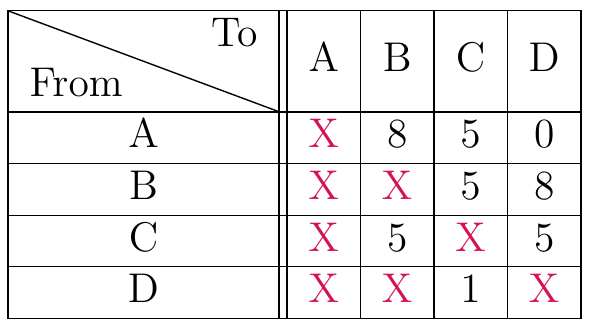}\\
    \caption{Example of a SOP instance: transition costs $c_{row,col}$, with \textcolor{safeA2}{X} indicating infeasible edges.}
    \label{fig:basicSOP}
\end{figure}
\begin{figure}[h!]
    \centering
    \begingroup
        \setlength{\tabcolsep}{20pt}
        \begin{tabular}{c c c}
            \includegraphics[scale=.45]{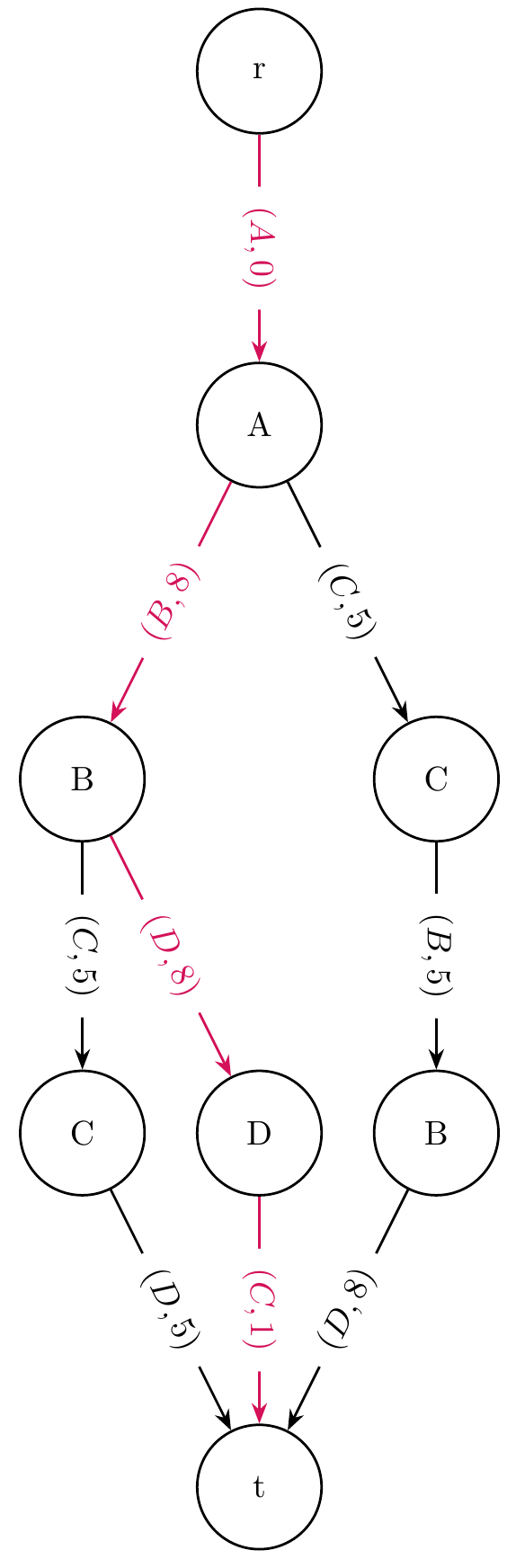} & \includegraphics[scale=.45]{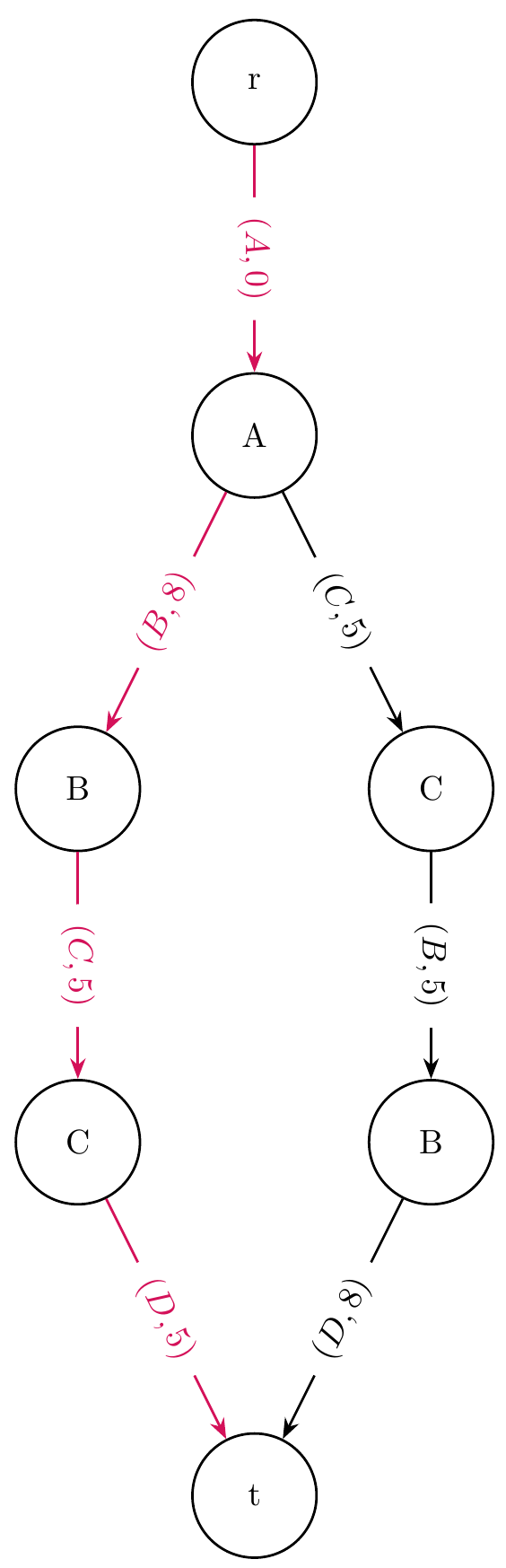} & \includegraphics[scale=.45]{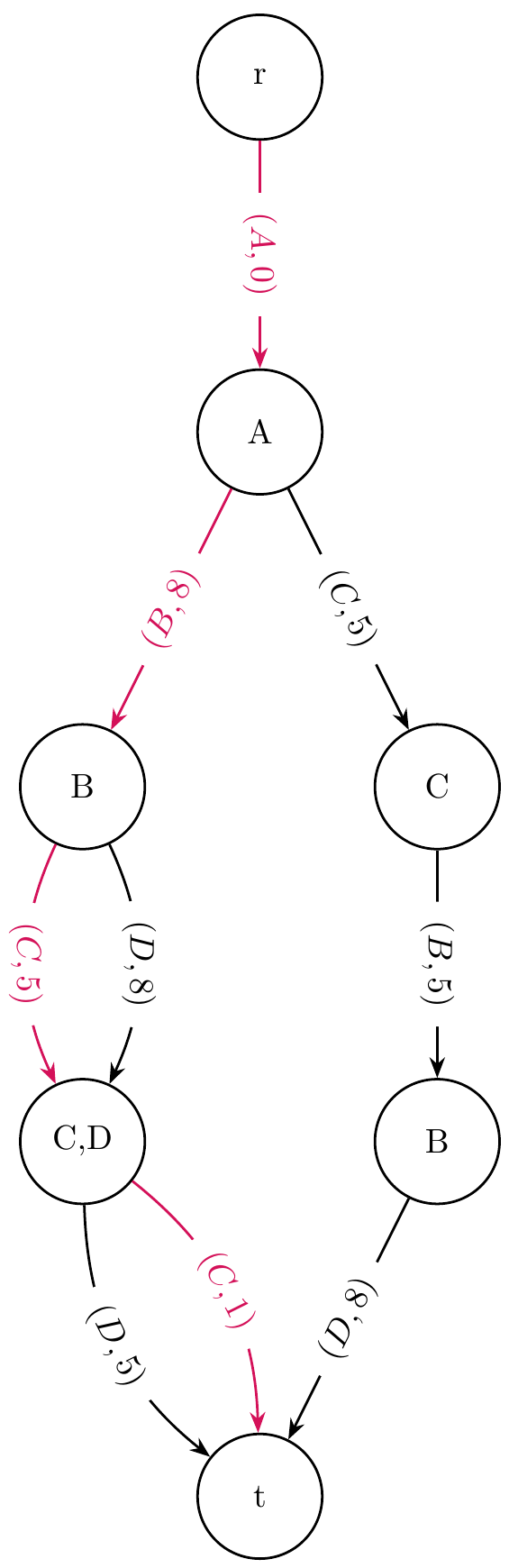}\\
            Exact Diagram & Restricted Diagram & Relaxed Diagram\\
            $z^* = [A,B, D,C]$ & $z^* = [A,B, C,D]$ & $z^* = [A,B, C,C]$\\
            $T^*=17$ & $T^*=18$ & $T^*=14$
        \end{tabular}
    \endgroup
    \caption{MDD Representation for the SOP instance presented in Figure \ref{fig:basicSOP}. Each arc $a$ has the format: $(l, v(a))$. The red path in each diagram indicates the shortest path from $r$ to $t$.}
    \label{fig:basicDDs}
\end{figure}

Figure \ref{fig:basicSOP} gives an instance of \emph{sequence ordering problem} (SOP), and Figure \ref{fig:basicDDs} contains simple examples of exact, restricted, and relaxed decision diagrams for that instance where $w_m=2$ for $\overline{\mathcal{M}}$ and $\underline{\mathcal{M}}$. The SOP requires finding the minimum-cost sequence of $n$ elements that includes each element exactly once, subject to transition costs $c_{ij}$ of following $x_i$ with $x_j$, and subject to precedence constraints requiring that certain elements precede others in the sequence. In other words, the SOP is an asymmetric TSP with precedence constraints. The label of each node matches the union of the labels of the incoming arcs. Each arc $a_{uv}$ is labeled in the format $(l, val(a_{uv}))$, representing the assignment of $x_{\ell_u}$ to $l$, and $val(a)$ represents the cost of the shortest path from the label of $u$ to $l$. In other words, an arc with label $l$ leaving layer $i$, represents the assignment of $l$ to the $i^{th}$ position of the sequence. The red path in each diagram indicates the shortest path through the diagram, and $T^*$ indicates the cost of the shortest path through the diagram.

\subsection{Restricted Decision Diagrams}
\label{sec:background:restricted}
Constructing $\overline{\mathcal{M}}$ for a given width $w_m$ is a straightforward process that can be thought of as a generalized greedy algorithm. Beginning with the root node $r$, an arc is generated for every element in the domain of $r$, and a node is generated at the end of each arc in the second layer. The process is repeated for each layer, except layer $n$ where all outgoing arcs point to the terminal, unless $w(\overline{\mathcal{M}})$ exceeds $w_m$. Then the least promising node is removed from the offending layer until $w(\overline{\mathcal{M}})$ is equal to $w_m$. The definition of \emph{least promising} is a heuristic decision. For the purposes of this paper, the least promising node is the node $u$ such that the shortest path from $r$ to $u$ is longer than the shortest path from $r$ to any other node $v\neq u$ in layer $\ell_u$.

It is of note that another method of reducing the width of $\overline{\mathcal{M}}$ is merging equivalent nodes. In the SOP, two nodes can be considered equivalent if they have the same state (last element in the sequence), and all incoming paths have visited the same set of elements. For example, a node with exactly one incoming path $[A,B,C]$ could be merged with a node in the same layer with exactly one incoming path $[B,A,C]$. In many MDD applications this is a valuable insight, and it helps motivate the algorithm for constructing $\underline{\mathcal{M}}$. However, for the SOP, we observed that the work of finding equivalent nodes in $\overline{\mathcal{M}}$ often outweighed the benefit of being able to merge nodes.

\subsection{Relaxed Decision Diagrams}
\label{sec:background:relaxed}

There are many methods of constructing relaxed decision diagrams, and many heuristic decisions that must be made when doing so. In this paper, we focus on the method described by Cire and van Hoeve~(2013)~\cite{MDDForSP} for sequencing problems. As opposed to the top-down construction described in Section \ref{sec:background:restricted}, here $\underline{\mathcal{M}}$ will be constructed by separation. Constructing DDs by separation uses $\underline{\mathcal{M}}$ as a domain store over which constraints can be propagated. This method starts with a weak relaxation, and then strengthens it by splitting nodes until each layer is either exact, or has a width equal to $w_m$. The algorithm begins from a $1$-width MDD with an arc from the node on layer $\ell_i$ to the node on layer $\ell_{i+1}$ for each element that can be placed at position $\ell_i$ in the sequence. Thus, even though each layer only has one node, there can be several arcs between layers (see the relaxed diagram in Figure \ref{fig:basicDDs}). Then a node $u$ is selected and split to strengthen the relaxation. The process of splitting $u$ involves creating two new nodes $u^\prime_1$ and $u^\prime_2$, and then distributing the \emph{in} arcs of $u$ between $u^\prime_1$ and $u^\prime_2$. Then for each \emph{out} arc $a_{uv}$ from $u$, arcs $a_{u^\prime_1v}$ and $a_{u^\prime_2v}$ are added such that $a_{uv}$, $a_{u^\prime_1v}$ and $a_{u^\prime_2v}$ all have the same label. Finally $u^\prime_1$ and $u^\prime_2$ are filtered to remove infeasible and sub-optimal arcs. A collection of filtering rules are used to check each arc. As an example, given a feasible solution to $\mathcal{P}$ with objective value $z_{opt}$, an arc $a$ can be removed if all paths containing $a$ have an objective value greater than $z_{opt}$. The full process of identifying which arcs can be removed is detailed in Cire and van Hoeve~(2013)~\cite{MDDForSP}, and is not replicated here.

The following notation and definitions are critical to understanding these algorithms. Let $All_u^\downarrow$ be the set of arc labels that appear in every path from $r$ to $u$. Let $Some_u^\downarrow$ be the set of arc labels that appear in at least one path from $r$ to $u$. Let $All_u^\uparrow$ and $Some_u^\uparrow$ be defined as above, except that they refer to paths from $u$ to $t$. Let $\mathscr{J}$ be the set of all possible arc labels. For the SOP, we define an \emph{exact} node $u$ as a node where $Some_u^\downarrow = All_u^\downarrow$ and all arcs ending at $u$ originate from exact nodes. Intuitively, a node $u$ is exact if all paths to $u$ contain the same set of labels, and all parents of $u$ are exact. Algorithm \ref{algo:RelaxedAlgo} formalizes the process of strengthening $\underline{\mathcal{M}}$.

\begin{algorithm}[h!t]
\SetAlgoLined
 Let $\underline{\mathcal{M}}$ be an MDD such that $Sol(\underline{\mathcal{M}}) \supseteq Sol(\mathcal{P})$\\
 \For{layer $L_j \in \underline{\mathcal{M}}$ from $j=1$ to $j=n$}{
 	 \While{$|L_j| < w_m$ and $\exists$ some node $ y \in L_j$ such that $y$ is not exact}{
 	 	$\mathscr{J} \leftarrow $ getAssignmentOrdering($\mathcal{P}$)\\
 	 	\small\emph{The getAssignmentOrdering() function returns a heuristically defined ordering of the values that can be assigned to decision variables}\normalsize\\
 	 	\For{$\phi \in \mathscr{J}$ \textbf{while} $|L_j| < w_m$}{
 			$S \leftarrow $ selectNodes($L_j$,$\phi$)\\
 			\small\emph{The selectNodes() function returns the set of nodes $u \in L_j$ such that $\phi \in Some_u^\downarrow \backslash All_u^\downarrow$}\normalsize\\
 			\For{$u \in S$ \textbf{while} $|L_j| < w_m$}{
 				Create two new nodes $u_1^\prime,u_2^\prime$\\
 				$L_j \leftarrow (L_j \cup \{u_1^\prime,u_2^\prime\})$\\
 				\ForEach{arc $a_{vu}$}{
					\eIf{$\phi \in (All_v^\downarrow ~ \cup$ the label of $a)$}
						{Redirect $a$ such that $a_{v u_1^\prime}$}
						{Redirect $a$ such that $a_{vu_2^\prime}$} 		
 				}
 				\ForEach{arc $a_{uv}$}{
 				    Create arcs $a_{u_1^\prime v}$ and $a_{u_2^\prime v}$ such that $label(a_{uv}) = label(a_{u_1^\prime v}) = label(a_{u_2^\prime v})$\\
 					filter($a_{u_1^\prime v}$), filter($a_{u_2^\prime v}$)\\
 					\small\emph{filter($a$) runs a list of quick checks to see if an arc can be removed}\normalsize\\
 				}
 				$L_j \leftarrow (L_j\backslash u)$\\
 			}
 		}
 	}
 }
 \Return{$\underline{\mathcal{M}}$}
 \caption{Refining Decision Diagrams for Sequencing~\cite{MDDForSP}}
  \label{algo:RelaxedAlgo}
\end{algorithm}

Deciding which nodes to split, and how to split them, are heuristic decisions with a significant impact on the bound that can be achieved without exceeding $w_m$~\cite{DDforO}. The algorithm discussed here selects nodes that can be split into equivalency classes, such that every path to the new node contains a certain label. Deciding which equivalency classes to produce first is another heuristic decision. The details of ordering the importance of the labels are specific to the problem being solved, and are not discussed here. However, it is important to note that the ordering for this implementation is static, and does not change between iterations. 

\subsection{Branch-and-Bound with Decision Diagrams}
\label{sec:background:bnb}

In a typical branch-and-bound algorithm, the branching takes place by splitting on the domain of the  variables. With decision diagrams, the branching takes place on the nodes themselves by selecting a set of exact nodes to represent the problem. The solver outlined by Bergman et al.~\cite{BnB2} defines an exact node as a node $u$ for which every path from $r$ to $u$ ends in an equivalent state. As mentioned above, we can be more specific when applying this to sequencing problems, and define an exact node $u$ as a node where $Some_u^\downarrow = All_u^\downarrow$ and all arcs ending at $u$ originate from exact nodes. An \emph{exact cutset} is defined as a set of exact nodes that contain every path from $r$ to $t$. Let $\underline{\mathcal{M}}(u)$ be a relaxed decision diagram with root $u$, and let $\overline{\mathcal{M}}(u)$ be a restricted decision diagram with root $u$. The branch-and-bound algorithm for MDDs proceeds by selecting an exact cutset of $\underline{\mathcal{M}}$, and using each node $u$ in the cutset as the root for a new restricted decision diagram $\overline{\mathcal{M}}(u)$ and relaxed decision diagram $\underline{\mathcal{M}}(u)$. A node can be removed from the queue if the relaxation of that node is not better than the best known solution to $\mathcal{P}$, otherwise the exact cutset of the new node is added to the queue, and the process repeats until the queue is empty. This is detailed by Algorithm \ref{algo:BnB}.

Gillard et al.~\cite{ddo} expanded on Algorithm \ref{algo:BnB} by incorporating a local search. A heuristic is used to quickly calculate a \emph{rough relaxed bound}\footnote{Gillard et al.~\cite{ddo} call the value \emph{rough upper bound}, but since we are testing a minimization problem in this paper, we use the term \emph{rough relaxed bound} instead.} at each node, and if the length of the shortest path to that node plus the rough relaxed bound is worse than the best known solution, the node can be removed. More formally, let $rrb(u)$ be a rough relaxed bound on $\mathcal{P}$ starting from node $u$, and let $z_{opt}$ be the value of best known solution so far. If $T^*(u)+rrb(u) > z_{opt}$, the node can be removed. They also provide evidence that if $rrb(u)$ is inexpensive to compute, it can be used to filter nodes in $\overline{\mathcal{M}}$ and $\underline{\mathcal{M}}$. The method of using a rough relaxed bound to trim nodes is used in this paper, but the details are problem specific and are discussed in a later section. 

\begin{algorithm}[!ht]
\SetAlgoLined
Let $\underline{\mathcal{M}_{uu^\prime}}$ be a partial diagram with root $u$ and terminal $u^\prime$\\
Let $v^*(u)$ be the lower bound of $\mathcal{P}$ resulting from starting at node $u$\\
Let $z_{opt}$ be the value of the best known solution\\
 $Q = \{r\}$ \\
 $v^*(r) \leftarrow 0$\\
 $z_{opt} \leftarrow \infty$\\
 \While{$Q \neq \emptyset$}{
 	$u \leftarrow$selectNode($Q$),  $Q \leftarrow Q\backslash\{u\}$\\
 	$\overline{\mathcal{M}} \leftarrow \overline{\mathcal{M}}(u)$\\
 	 \If{$v^*(\overline{\mathcal{M}}) < z_{opt}$} {$z_{opt} \leftarrow v^*(\overline{\mathcal{M}})$\\}
 	 \If{$\overline{\mathcal{M}}$ is not exact}{
 	 	$\underline{\mathcal{M}} \leftarrow \underline{\mathcal{M}}(u)$\\
 	 	\If{$v^*(\underline{\mathcal{M}}) < z_{opt}$} {
 	 		$S \leftarrow $ exactCutset($\underline{\mathcal{M}}$)\\
 	 		\ForEach{$u^\prime \in S$}{
				let $v^*(u^\prime) = v^*(u) + v^*(\underline{\mathcal{M}_{uu^\prime}})$\\
				$Q \leftarrow Q \cup u^\prime$\\
 	 		}
 	 	}
 	 }
 }
 \Return{$z_{opt}$}
 \caption{Decision Diagram based Branch-and-Bound  (BnB)~\cite{BnB2}}
 \label{algo:BnB}
\end{algorithm}

\section{Peel-and-Bound Algorithm}
\label{sec:pnb}

The motivation for peel-and-bound stems from an observation about Algorithm~\ref{algo:RelaxedAlgo}. When implemented in a branch-and-bound structure, a large portion of the work done while generating each $\underline{\mathcal{M}}$ is repeated at every iteration. Creating the relaxation for some exact node $u$ in the queue requires creating a $1$-width decision diagram, iterating over each layer from the top down, and splitting nodes in a predetermined order. The static order of node splits means that for each node $y$ such that $\ell_y > \ell_u$, the first equivalency class created when splitting $y$ is the same in $\underline{\mathcal{M}}(r)$ and $\underline{\mathcal{M}}(u)$. The existing arcs for both diagrams will be sorted in the same way, and the only difference is the possibility of filtering arcs in $\underline{\mathcal{M}}(u)$ that could not be filtered in $\underline{\mathcal{M}}(r)$ due to the added constraint that all paths must pass through $u$. The extra filtered arcs are the reason that $\underline{\mathcal{M}}(u)$ may produce a stronger bound than $\underline{\mathcal{M}}(r)$. However, because equivalency classes are chosen in the same order each time, many arcs that were filtered while constructing $\underline{\mathcal{M}}(r)$ will also be filtered again while constructing $\underline{\mathcal{M}}(u)$. There is a sub-graph of $\underline{\mathcal{M}}(r)$, induced by node $u$, that contains all of the paths that will be encoded in $\underline{\mathcal{M}}(u)$, but does not contain the arcs that are filtered from both diagrams during construction. Thus, less work needs to be performed at each iteration of branch-and-bound by starting from that sub-graph instead of a $1$-width diagram. If the split order is static, the same diagram is generated starting from either the $1$-width diagram, or the sub-graph induced by $u$. If the split order changes between branch-and-bound iterations, the sub-graph induced by $u$ is still a valid relaxation, but the generated diagram will differ from one that began at width $1$. 

\begin{figure}[!ht]
    \centering
    \begingroup
        \setlength{\tabcolsep}{10pt}
        \begin{tabular}{c c c c}
            \includegraphics[scale=.4]{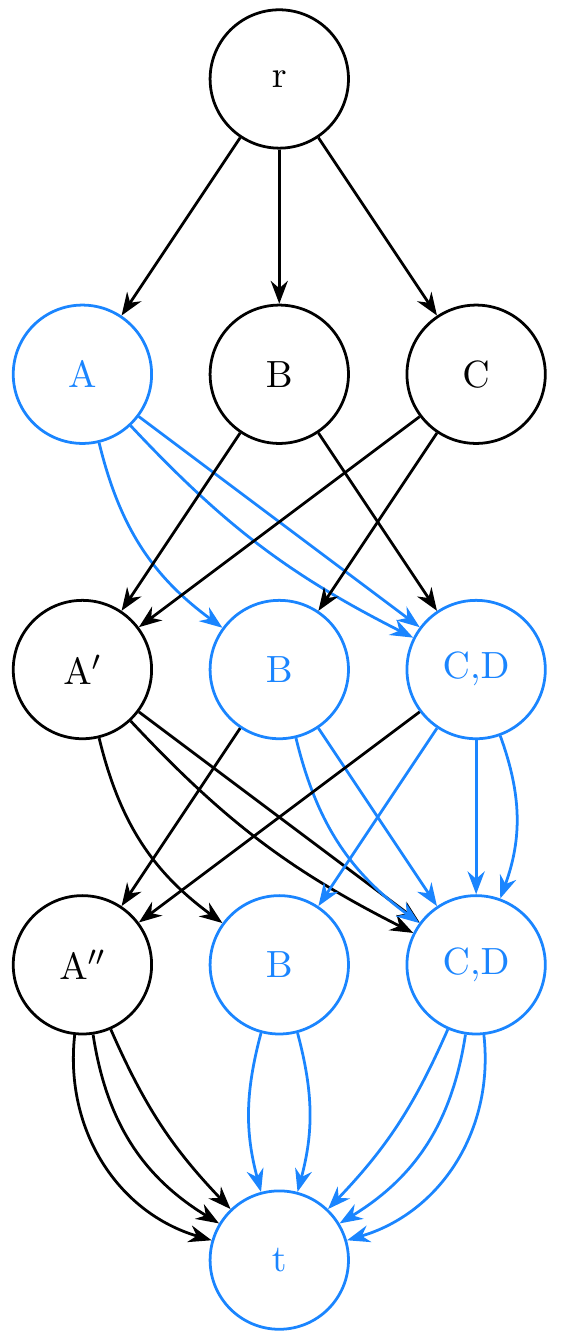} & \includegraphics[scale=.4]{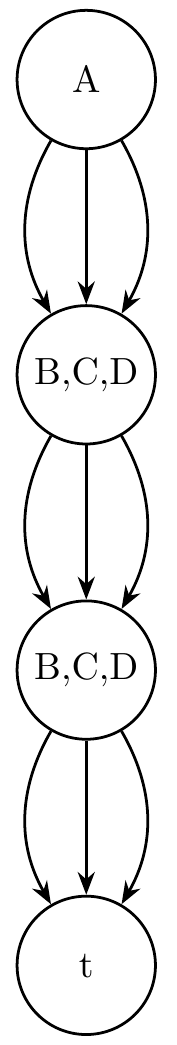} & \includegraphics[scale=.4]{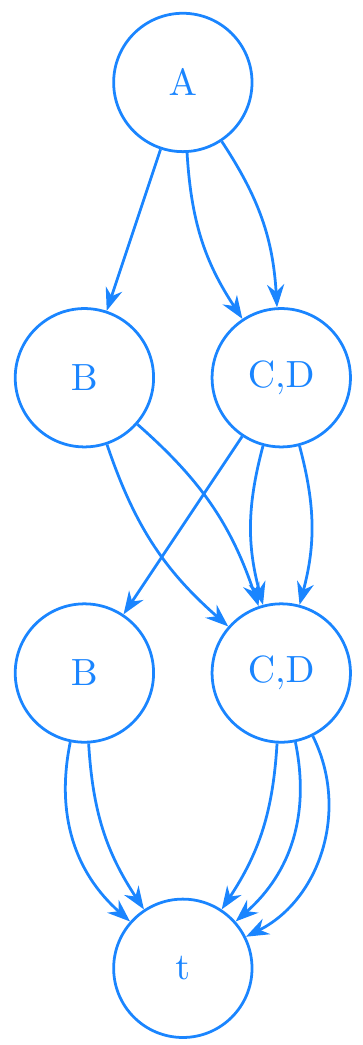} & \includegraphics[scale=.4]{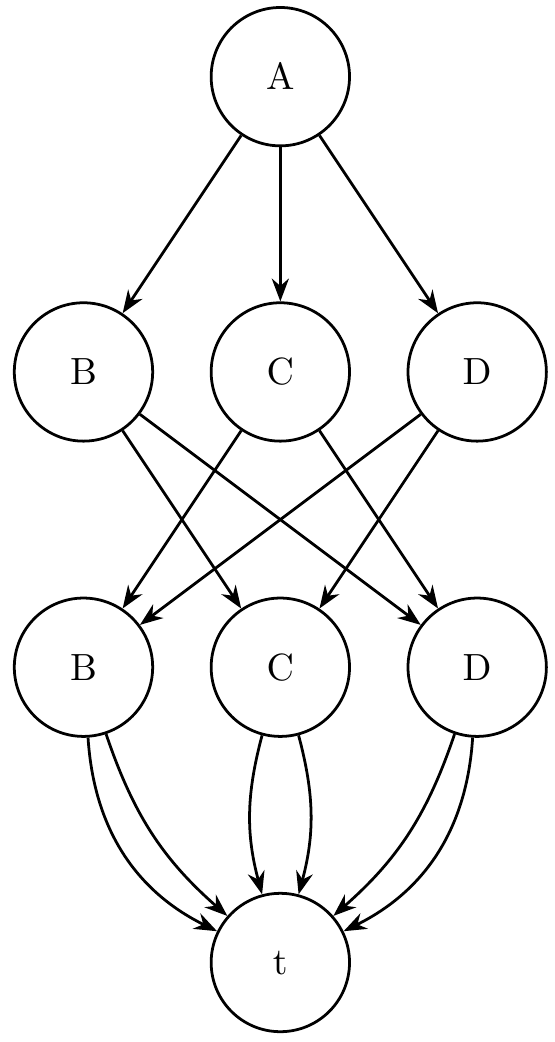}\\
            $\underline{\mathcal{M}}(r)$ & $\underline{\mathcal{M}}(A)$ Stage 1 & $\underline{\mathcal{M}}(A)$ Stage 2 & $\underline{\mathcal{M}}(A)$ Stage 3
        \end{tabular}
    \endgroup
    \caption{Example of an induced sub-graph for a SOP instance (shown in blue), and the associated relaxed decision diagram with the same root.}
    \label{fig:SubDiagram}
\end{figure}

Consider a SOP instance where the goal is to order the elements $[A,B,C,D]$, subject to the precedence constraint that $A$ must precede $D$, an alphabetical ordering heuristic, and $w_m = 3$. Figure~\ref{fig:SubDiagram} shows $\underline{\mathcal{M}}(r)$, and $\underline{\mathcal{M}}(A)$ in three stages. The first stage is the initial $1$-width diagram. The second stage is after one split on each layer, and the third stage is the complete diagram. The sub-graph shared by $\underline{\mathcal{M}}(r)$ and $\underline{\mathcal{M}}(A)$ is highlighted in blue, indicating that in this case the first two splits could have been read from $\underline{\mathcal{M}}(r)$ instead of being re-created from scratch. For the sake of legibility, arc values and arc labels are not included.

\begin{figure}[!ht]
    \centering
    \begingroup
        \setlength{\tabcolsep}{1pt}
        \begin{tabular}{c c c c}
            \includegraphics[scale=.36]{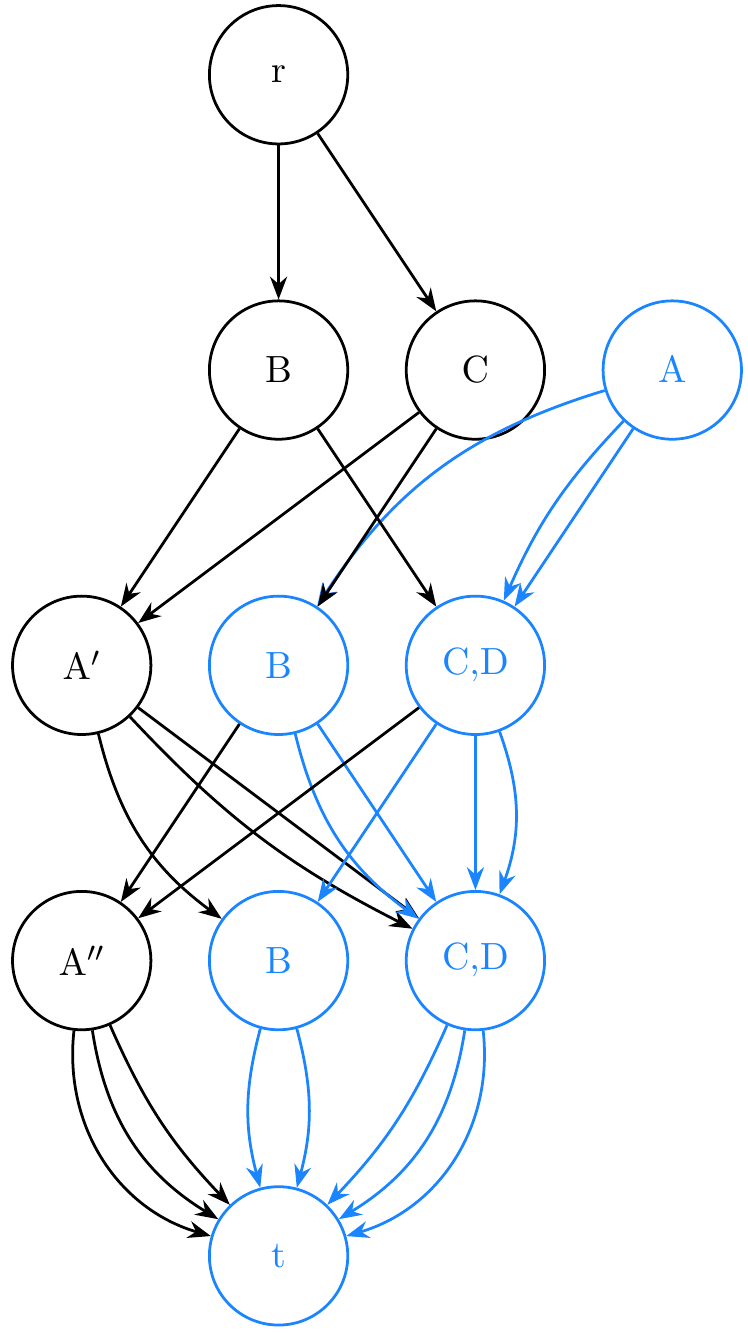} & \includegraphics[scale=.36]{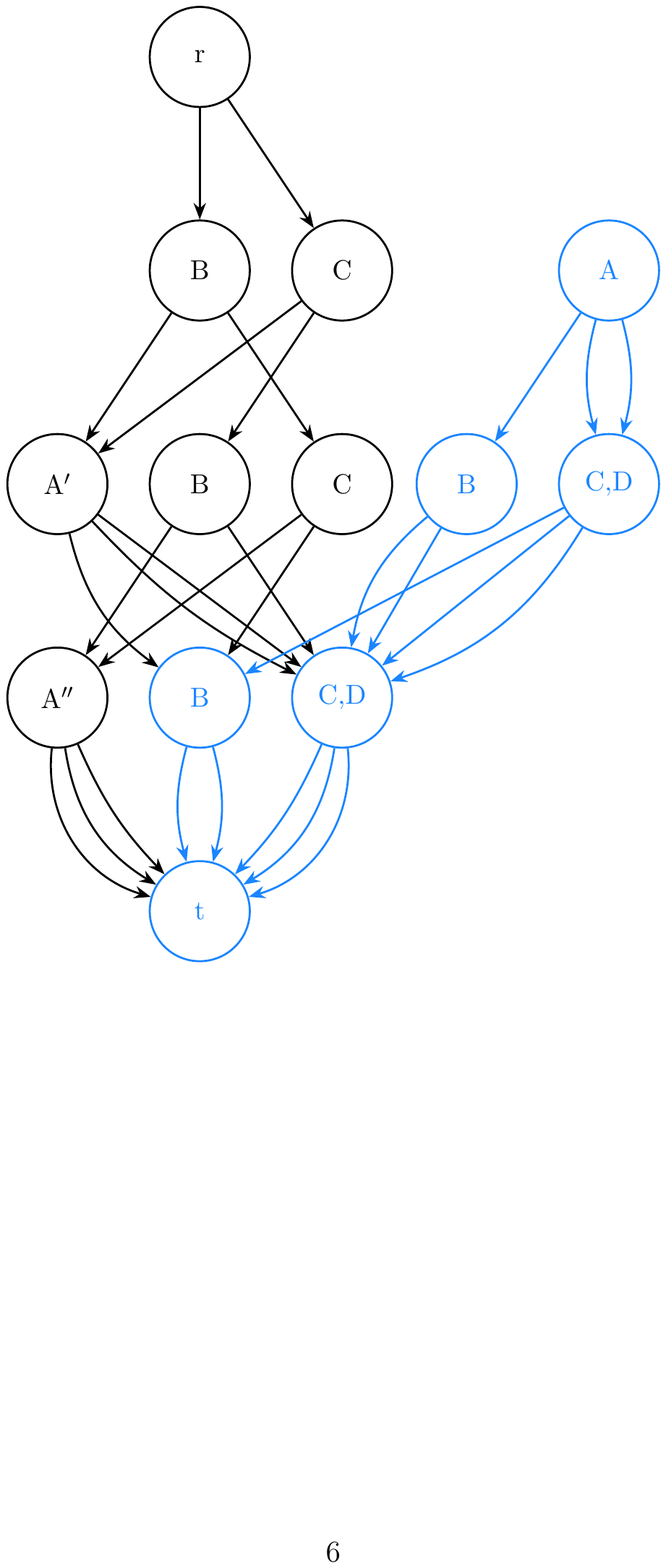} & \includegraphics[scale=.36]{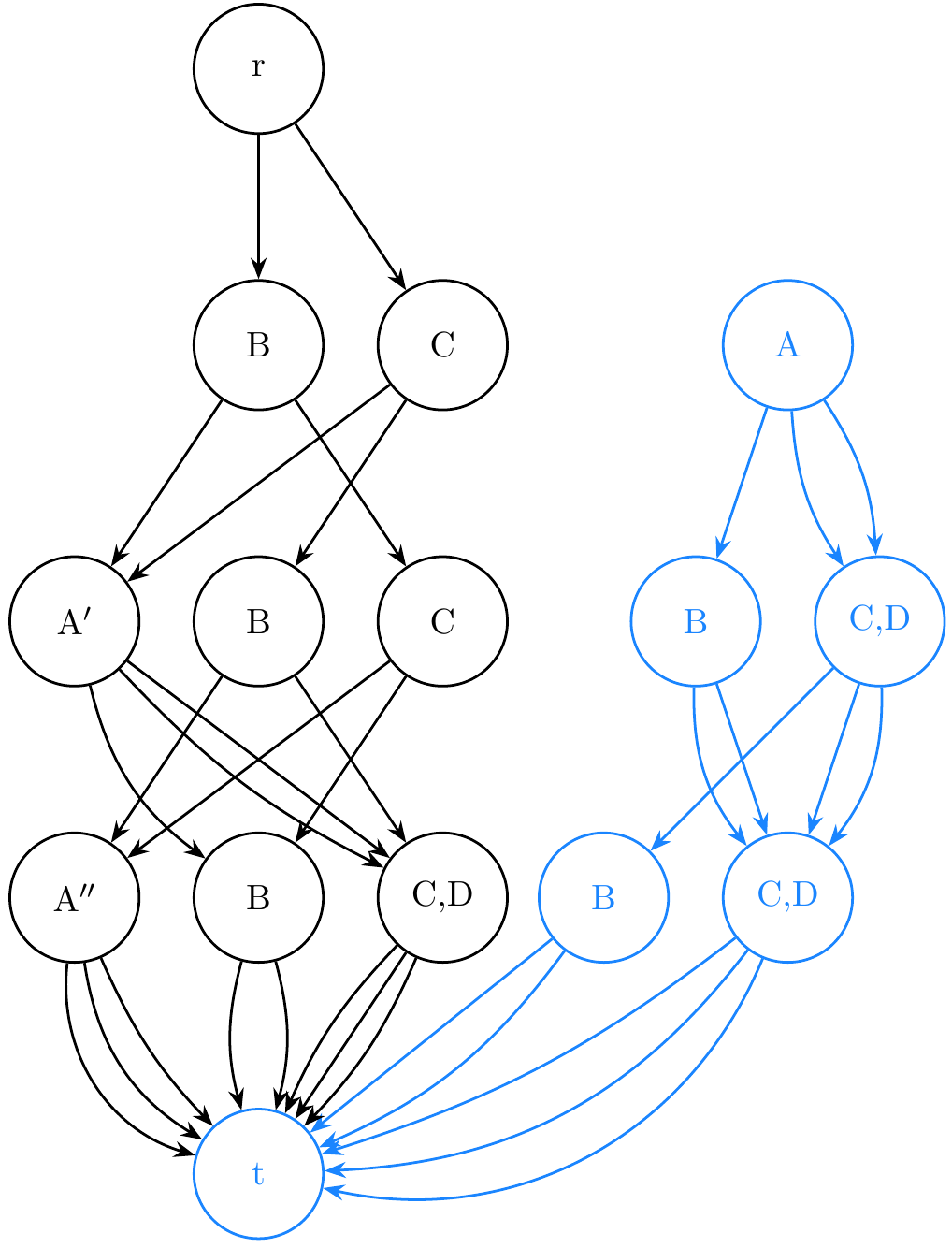} & \includegraphics[scale=.36]{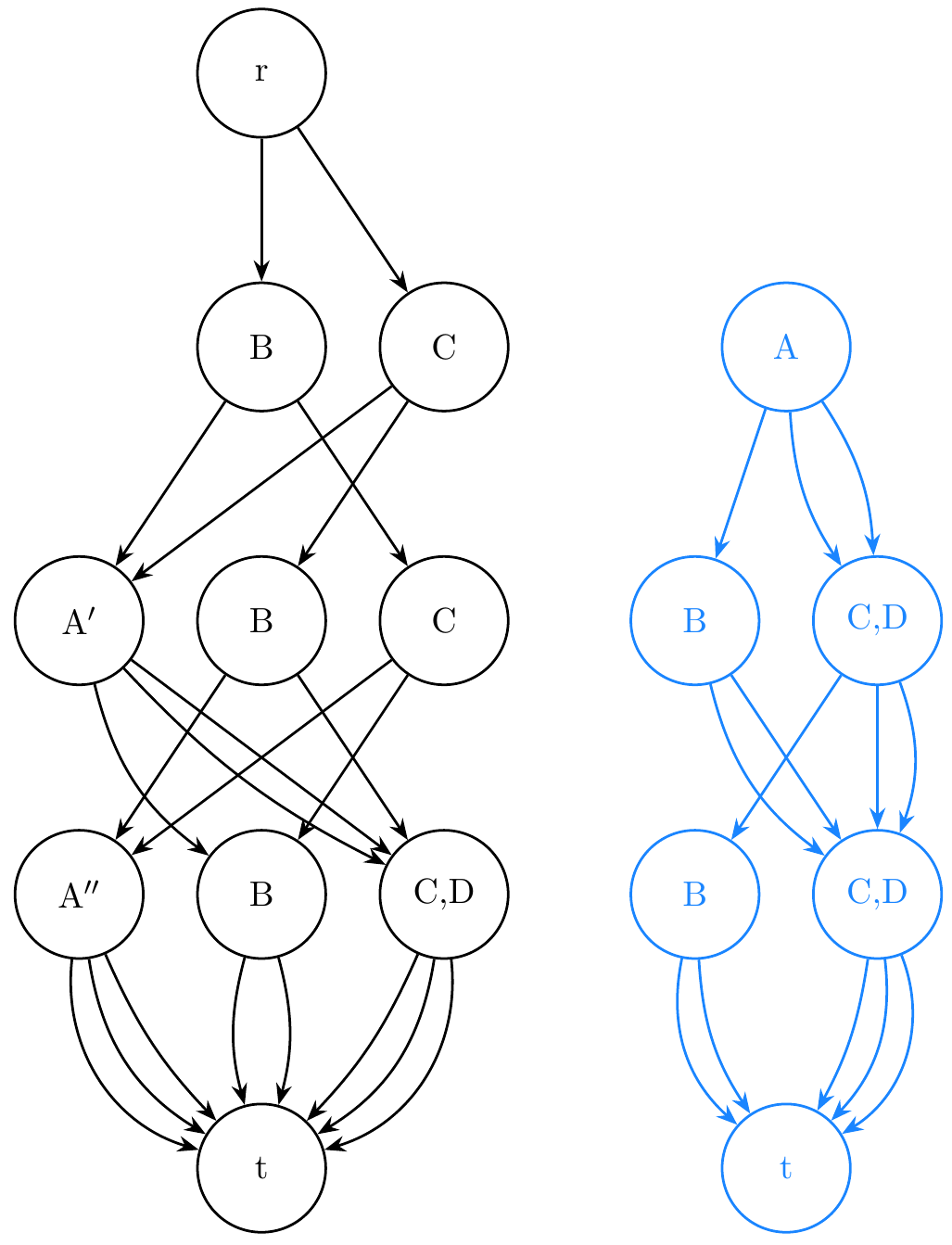} \\
            (1) & (2) & (3) & (4)\\
        \end{tabular}
    \endgroup
    \caption{An example of a peel operation. In (1), $A$ is selected to induce the peel process and removed from the the original diagram ($\underline{\mathcal{M}}(r)$ from Figure \ref{fig:SubDiagram}). In (2) the arcs that connect $A$ to the original diagram are moved to copies of the nodes they originally ended at, and infeasible arcs are filtered. In (3) and (4) the process is repeated until the diagrams are disconnected.}
    \label{fig:PnB}
\end{figure}

This mechanism can be embedded into a slightly modified version of the standard branch-and-bound algorithm based on decision diagrams (Algorithm~\ref{algo:BnB}). In peel-and-bound, the queue stores diagrams instead of nodes. After the initial relaxation $\underline{\mathcal{M}}(r)$ is generated, the entire diagram is placed into the queue $Q$ such that $Q = \{\underline{\mathcal{M}}(r)\}$. Then a diagram $\underline{\mathcal{M}}(u)$ is selected from $Q$ (for the first iteration $\underline{\mathcal{M}}(u) = \underline{\mathcal{M}}(r)$). However, instead of selecting an exact cutset of $\underline{\mathcal{M}}(u)$, a single exact node $e$ from $\underline{\mathcal{M}}(u)$ is selected. The process of selecting a diagram and exact node are heuristic decisions that are discussed in Section~\ref{sec:pnb:implementation}. The process of \emph{peeling} $e$ is as follows. Create an empty diagram $\underline{u}$, remove $e$ from $\underline{\mathcal{M}}(u)$, and then put $e$ into $\underline{u}$ such that $e$ is the root of $\underline{u}$, and the arcs leaving $e$ still end in $\underline{\mathcal{M}}(u)$. Then for each node $y$ in $\underline{\mathcal{M}}(u)$ with an \emph{in} arc that originates in $\underline{u}$, a new node $y^\prime$ is made and added to $\underline{u}$. Each \emph{in} arc $a_{oy}$ of $y$ that originates in $\underline{u}$ is removed and then arc $a_{oy^\prime}$ is added to $\underline{u}$. Then the \emph{out} arcs of $y$ and $y^\prime$ are filtered using the same \emph{filter} function as Algorithm~\ref{algo:RelaxedAlgo}. The process of removing and adding arcs is repeated until there are no arcs ending in $\underline{\mathcal{M}}(u)$ that originate in $\underline{u}$. This procedure accomplishes a top-down reading of the sub-graph induced by $e$, and potentially strengthens $\underline{\mathcal{M}}(u)$ by removing nodes and arcs in the process. If the shortest path through the modified $\underline{\mathcal{M}}(u)$ is less than the best known solution, $\underline{\mathcal{M}}(u)$ is put back into $Q$. $\underline{u}$ is relaxed using Algorithm \ref{algo:RelaxedAlgo}, let $\underline{\mathcal{M}}(\underline{u})$ be the result; then if the shortest path through the refined diagram $\underline{\mathcal{M}}(\underline{u})$ is less than the best known solution, $\underline{\mathcal{M}}(\underline{u})$ is added to $Q$. The whole procedure is repeated until there are no nodes left in the queue ($Q = \emptyset$). A peel operation is illustrated and explained in Figure \ref{fig:PnB}. Peel-and-bound is formalized in 
Algorithm~\ref{algo:PnB}, and the peel operation is formalized in Algorithm~\ref{algo:Peel}.

\begin{algorithm}[!ht]
\SetAlgoLined
Let $v^*(u)$ be the lower bound of $\mathcal{P}$ resulting from starting at node $u$\\
Let $z_{opt}$ be the value of the best known solution\\
 $Q = \{\underline{\mathcal{M}}(r)\}$ \\
 $z_{opt} \leftarrow \infty$\\
 \While{$Q \neq \emptyset$}{
    $\mathscr{D} \leftarrow$selectDiagram($Q$), $Q \leftarrow Q\backslash\{\mathscr{D}\}$\\
    $u \leftarrow$selectExactNode($\mathscr{D}$)\\
    $\underline{u},\mathscr{D}^*\leftarrow$ peel($\mathscr{D}$, $u$) \emph{(See Algorithm \ref{algo:Peel})}\\
    \If{$v^*(\mathscr{D}^*) < z_{opt}$}{
        $Q \leftarrow Q \cup \{\mathscr{D}^*\}$
    }
    $\overline{\mathcal{M}} \leftarrow \overline{\mathcal{M}}(u)$\\
    \If{$v^*(\overline{\mathcal{M}}) < z_{opt}$} {$z_{opt} \leftarrow v^*(\overline{\mathcal{M}})$\\}
 	 \If{$\overline{\mathcal{M}}$ is not exact}{
 	 	$\underline{\mathcal{M}} \leftarrow \underline{\mathcal{M}}(\underline{u})$\\
 	 	\If{$v^*(\underline{\mathcal{M}}) < z_{opt}$} {
 	 		$Q \leftarrow Q \cup \{\underline{\mathcal{M}}\}$
 	 	}
 	 }
 }
 \Return{$z_{opt}$}
 \caption{Peel-and-Bound (PnB) Algorithm}
 \label{algo:PnB}
\end{algorithm}

Separating each node $u$ during a peel requires creating a new node $u^\prime$, moving the \emph{in} arcs of $u$ that originate in the peeled diagram $\underline{u}$ to \emph{$u^\prime$}, copying the \emph{out} arcs of $u$ to $u^\prime$, and then filtering the \emph{out} arcs of $u$ and $u^\prime$. Creating a new node in our implementation has a time in $\mathcal{O}(n)$ due to storing state information that has a size in $\mathcal{O}(n)$ (such as $All_{u^\prime}^\downarrow$). However, it is possible that in other applications the size of a node is in $\mathcal{O}(1)$. The number of \emph{in} arcs of $u$ is at most $w$, although this worst case is unlikely in practice because it requires $\underline{u}$ to have width $w$ and for each node in $\underline{u}$ on layer $\ell_{u}-1$ to have an arc ending at $u$. Thus, moving the \emph{in} arcs of $u$ has a time in $\mathcal{O}(w)$. The number of \emph{out} arcs of $u$ is at most $n$, and each arc has a size in $\mathcal{O}(1)$, so copying the \emph{out} arcs has a time in $\mathcal{O}(n)$. Each individual filtering process has a time in $\mathcal{O}(1)$ as it uses only existing state information from $u$ and $u^\prime$, and it is performed on the at most $2n$ \emph{out} arcs of $u$ and $u^\prime$. Thus, filtering the \emph{out} arcs has a time in $\mathcal{O}(n)$. Therefore, separating one node during the peel process has a time in $\mathcal{O}(n+w)$. Separations during a standard relaxation procedure require selecting a node ($\mathcal{O}(w)$), making a new node ($\mathcal{O}(n)$), partitioning the \emph{in} arcs ($\mathcal{O}(nw)$), copying the \emph{out} arcs ($\mathcal{O}(n)$), and filtering the \emph{out} arcs ($\mathcal{O}(n)$). The reason that there can be more \emph{in} arcs during a standard relaxation procedure is because the nodes in a $1$-width diagram can have \emph{in} arcs with different labels coming from the same node, whereas the structure of the diagram during a peel guarantees that each node $u$ can have only one \emph{in} arc from each node on the layer $\ell_{u}-1$. Thus, the total time for a separation in a standard relaxation is in $\mathcal{O}(nw)$.

\begin{algorithm}[!ht]
\SetAlgoLined
    Let $in(u)$ for some node $u$ be the set of arcs that end at node $u$\\
    Let $out(u)$ for some node $u$ be the set of arcs that originate from node $u$\\
    Let $in(\mathcal{M})$ for some MDD $\mathcal{M}$ be the set of arcs that end in $\mathcal{M}$\\
    Let $out(\mathcal{M})$ for some MDD $\mathcal{M}$ be the set of arcs that originate in $\mathcal{M}$\\
    \textbf{input:} a relaxed MDD $\mathscr{D}$, and an exact node $u$ in $\mathscr{D}$\\
    Let $\underline{u}$ be an empty decision diagram\\
    $in(u) \leftarrow \emptyset$\\
    $\mathscr{D} \leftarrow \mathscr{D}\backslash u$\\
    $\underline{u} \leftarrow u$\\
    \While{$in(\mathscr{D})\cap out(\underline{u}) \neq \emptyset$}{
        \ForEach{node $m \in \mathscr{D}$ with an in arc that originates in $\underline{u}$}{
            create a new node $m^\prime$, and add it to $\underline{u}$\\
            \ForEach{arc $a \in in(m)$ that originates in $\underline{u}$}{
                change the destination of $a$ to $m^\prime$\\
                filter($a$)
            } 
            \ForEach{arc $a \in out(m)$}{
                filter($a$)
            }
        }
    }
    \While{$\exists$ some node $m \in D$ with $in(m)=\emptyset$ or $out(m) = \emptyset$ (excluding $r$ and $t$)}{
            $in(m) \leftarrow \emptyset$\\
            $out(m) \leftarrow \emptyset$\\
            $\mathscr{D} \leftarrow \mathscr{D} \backslash \{m\}$
    }
    \Return{($\underline{u}$, $\mathscr{D}$)}
    \caption{Peeling process used in Algorithm \ref{algo:PnB}}
    \label{algo:Peel}
\end{algorithm}

The maximum number of separations during a peel is the maximum number of nodes in the peeled diagram. A peeled diagram can have at most $(n-3) \times w + 2$ nodes, and thus the number of nodes is in $\mathcal{O}(nw)$. Therefore, the entire peel process has a time in $\mathcal{O}(n^2w + nw^2)$. The maximum number of separations during a standard relaxation is the exact same as during a peel, since the resulting diagram will be the same size. Thus, the standard relaxation has a total time in $\mathcal{O}(n^2w^2)$. However, peel-and-bound uses a peel to generate some fraction of the nodes, then a standard relaxation to generate the rest. Let $\alpha$ be the percent of nodes that are peeled during the peel. It follows that the total time for an iteration of peel-and-bound is in $\mathcal{O}(\alpha(n^2w + nw^2)+(1-\alpha)(n^2w^2))$. Therefore, the larger that $\alpha$ grows, the more time peel-and-bound saves over branch-and-bound.

\subsection{Advantages and Implementation Decisions}
\label{sec:pnb:implementation}
The branch-and-bound algorithm proposed by Bergman et al.~(2016)~\cite{BnB2} requires selecting an exact cutset of $\underline{\mathcal{M}}$. Peel-and-bound requires selecting a diagram from the queue, and an exact node to start the peel process. The choice of node has a substantial impact on how quickly the process converges to an optimal solution, because it serves two purposes simultaneously. As discussed earlier, the first purpose of peeling is to avoid recreating a portion of the diagram at each iteration. The second purpose is to strengthen the overall relaxation. Let $\underline{u}$ be a diagram peeled from $\underline{\mathcal{M}}$, and let $\underline{\mathcal{M}}^*$ be $\underline{\mathcal{M}}$ after the peel operation. If $Sol(\mathcal{P}) \subseteq Sol(\underline{\mathcal{M}})$ then $Sol(\mathcal{P}) \subseteq Sol(\underline{\mathcal{M}}^*)\cup Sol(\underline{u})$. The only step of peel-and-bound that removes paths is the \emph{filter} step, which only removes an arc if no feasible solutions can pass through that arc. If the node the peel is induced from contains the shortest path through $\underline{\mathcal{M}}$, then there will be a new shortest path through $\underline{\mathcal{M}}^*$ with $T^*(\underline{\mathcal{M}}^*) \geq T^*(\underline{\mathcal{M}})$. Similarly after peeling, the peeled diagram is going to be strengthened and $T^*(\mathcal{M}(\underline{u})) \geq T^*(\underline{u})$. Therefore, when implementing the \emph{selectDiagram} and \emph{selectExactNode} functions from Algorithm \ref{algo:PnB}, we propose selecting the diagram $\mathscr{D}$ with the weakest bound, and an exact node from $\mathscr{D}$ that contains $z^*(\mathscr{D})$ at each iteration. Using these parameters, the peel step of peel-and-bound strengthens the relaxed bound of $\mathcal{P}$, in addition to providing a stronger initial diagram to use when generating $\mathcal{M}(\underline{u})$.

We propose two heuristics for selecting a node from $\mathscr{D}$ that contains $z^*(\mathscr{D})$. The first heuristic picks the first node in the shortest path through the diagram with at least one child that is not exact, we call this the \emph{last exact node}. The second heuristic picks the \emph{frontier} node, the highest-index exact node that contains $z^*(\mathscr{D})$. Taking the last exact node is more of a breadth-first search, taking the largest possible set of nodes that can be strengthened (anything above the last exact node is exact, and cannot be improved). In contrast, taking the frontier node is more of a depth-first search, taking fewer nodes and exploring those nodes at greater depth.

Cire and van Hoeve~(2013)~\cite{MDDForSP} propose that each iteration of Algorithm \ref{algo:RelaxedAlgo} starts from a $1$-width MDD. However, for peel-and-bound with a non-separable objective function, starting from a $1$-width MDD poses a problem. The arcs in such a diagram do not have exact values, because they are dependent on the state of the node they originate from. As nodes are peeled, the values of those arcs must be updated, and the operation becomes computationally expensive at scale. This problem can be avoided by creating the initial diagram using a structure where all of the arcs ending at a given node have the same label. The resulting initial diagram has a width of $n$, and each node on the layer is assigned to one state $s \in \{1,...,n\}$. Then every possible feasible arc between consecutive layers is added. Thus, the nodes of $\underline{\mathcal{M}}$ do not have relaxed states, and each arc can only take one possible value. Starting from such a diagram not only removes the need to update arc values, it ensures that every arc generated during peel-and-bound is an exact copy of an arc that exists in the initial diagram, since arcs are only copied or removed, never updated or added. An alternative method of handling non-separable objective functions is explored by Hooker~\cite{canonicalarcs,jobsequencing,jobsequencing2}.

\subsection{Limitations and Handling Memory}
\label{sec:pnb:limitations}
The focus of this paper is sequencing problems, but peel-and-bound can be easily applied to other optimization problems. However, some existing MDD based methods conflict with peel-and-bound. For example, some MDD algorithms use a dynamic variable order~\cite{karahalios}, such that the variables the layers on $\overline{M}$ are mapped to in one iteration of branch-and-bound, are different in the next. Peel-and-bound as it is presented in this paper cannot be paired with a dynamic variable order. Furthermore, the method in this paper is specific to decision diagrams generated using separation. We believe the method can be extended to decision diagrams that use a merge operator, but it has not been shown here.

Memory limitations present a problem for peel-and-bound in theory, but not in practice. Each open diagram remains in the queue, and thus must be stored in memory. However, this problem can be handled in many ways; two are given here. A dynamic method of handling the problem is to start targeting large diagrams with bounds close to $z_{opt}$ as memory limitations start to become a problem. Such diagrams can usually be closed quickly, and subsequently removed from memory, freeing up space for the algorithm to continue. Alternatively, the diagrams with bounds closest to $z_{opt}$ can be deleted in favor of storing just the root node, then when they need to be processed, initial diagrams are generated for those once again. This method essentially falls back to branch-and-bound until memory limitations cease to be a problem. Additional approaches for working with memory limitations, and evidence that the problem can be handled efficiently, are presented in Perez and Régin~(2018)~\cite{Perez_Regin_2018}.

\subsection{Integrating Rough Relaxed Bounds}
\label{sec:pnb:rrb}
This implementation incorporates the rough relaxed bounding method proposed by Gillard et al.~\cite{ddo}. Rough relaxed bounding was used to trim the domain of each node during construction of the restricted DDs, and was also added as a check to the \emph{filter} function in Algorithm~\ref{algo:RelaxedAlgo}. When the initial model is created, a map is also created from each node $u$, to a list of the other nodes sorted by their distance from $u$. The rough relaxed bound $rrb(a)$ of an arc $a_{fg}$ was calculated as follows. For each node $u$ that has not necessarily been visited ($u \notin All_g^\downarrow$), look up the shortest distance from that node to a different node that has also not been visited. Then,  sort the resulting list, and repeatedly remove the largest value until the list has a length equal to the number of remaining decisions. The sum of the values in the list, plus the value of the shortest path from $r$ to the end of $a$, is the rough relaxed bound of $a$. If $rrb(a)$ is worse than the best known solution, the arc is removed.

\section{Experiments on the Sequence Ordering Problem}
\label{sec:results}
The goal of this section is to assess the performances of the peel-and-bound algorithm (PnB, Algorithm~\ref{algo:PnB}). 
To do so, a comparison with the standard decision diagram based branch-and-bound algorithm (BnB, Algorithm~\ref{algo:BnB}) is proposed. Both algorithms are implemented in Julia and are open-source\footnote{https://github.com/IsaacRudich/PnB$\_$SOP}. 
To ensure a fair comparison, both algorithms resort to the same function for generating relaxed decision diagrams (Algorithm~\ref{algo:RelaxedAlgo}), and the same function for generating restricted decision diagrams. While the functions being called are the same, there are two differences at run-time. At the end of line $26$ in Algorithm~\ref{algo:RelaxedAlgo}, an additional operation runs during BnB where the values of the arcs leaving layer $j$ are updated. The second difference is that BnB starts each relaxation from a $1$-width DD, while PnB passes a partially completed diagram to the relaxation function as a starting point.

The testing environment was built from scratch to ensure a fair comparison, 
so it lacks the many propagators used by cutting-edge solvers like CPO to remove nodes from the PnB/BnB queue~\cite{cpS, MDDForSP}. However, it provides a clean comparison of the two algorithms by requiring that every function used by both BnB and PnB is exactly the same between the two, with the only differences arising due to PnB's ability to ensure that all arcs are exact from the beginning. All of the heuristic decisions that were made are identical for both algorithms. 

\subsection{Description of the Heuristics Considered}
\label{sec:results:heuristics}
The \emph{sequence ordering problem} can be considered as an asymmetric \emph{travelling salesperson problem} with precedence constraints. The objective is to find a minimum cost path that visits each of the $n$ elements exactly once, and respects the precedence constraints. The method used for generating relaxed DDs requires creating a heuristic ordering of all possible arc assignments by importance. The arc values in this case are representative of the $n$ elements in the path. The ordering used was generated by sorting the $n$ elements, first by their average distance from the other elements, and then by the number of elements each element must precede. The resulting order places a higher importance on elements that are far away from other elements and must precede many other elements.

The branch-and-bound algorithm processes nodes in an order designed to try and improve the existing relaxed bound at each iteration. When a node $u$ is added to the BnB queue, it is assigned a value equal to the value of the shortest path from the root $r$ to the terminal $t$, that passes through $u$. The best known relaxed bound on the problem is the smallest value of a node in the queue, and that node is always chosen to be processed. Peel-and-bound is implemented with the same goal of improving bounds at each iteration. However, PnB stores diagrams, not nodes. Let the value of a diagram be the value of the shortest path to the terminal. At each iteration of peel-and-bound, the diagram with the lowest value is selected, and then a node is chosen from that diagram to induce the peel process. All of the experiments here used a process where the selected node is the first node in the shortest path from $r$ to $t$ with a child node that is not exact (the last exact node). Testing was done to determine whether using the last exact node or the frontier node would perform better for the problem being considered, but there was not a significant difference between the two during any of the tests. Several of the benchmark problems were run using various decision diagram widths, and the last exact node was chosen because it sometimes showed a very slight improvement over the frontier node. While it is likely that this choice makes a difference on some problems, it does not matter for the SOP. 

\subsection{Experimental Results}
\label{sec:results:experiment}
The experiments were performed on a computer equipped with an AMD Rome 7532 at 2.40 GHz with 64Gb RAM. The solver was tested using DD widths of $64, 128$, and $256$ on the $41$ SOP problems available in TSPLIB~\cite{tsplib}. For comparisons between PnB and BnB, a timestamp, new bounds, and the length of the remaining queue were recorded each time the bounds on a problem were improved. Another experiment was performed to test the scalability of PnB at width $2048$, for which only the final bounds were recorded. Execution time was limited to $3,600$ seconds.

\begin{figure}[!ht]
    \includegraphics[scale=.53]{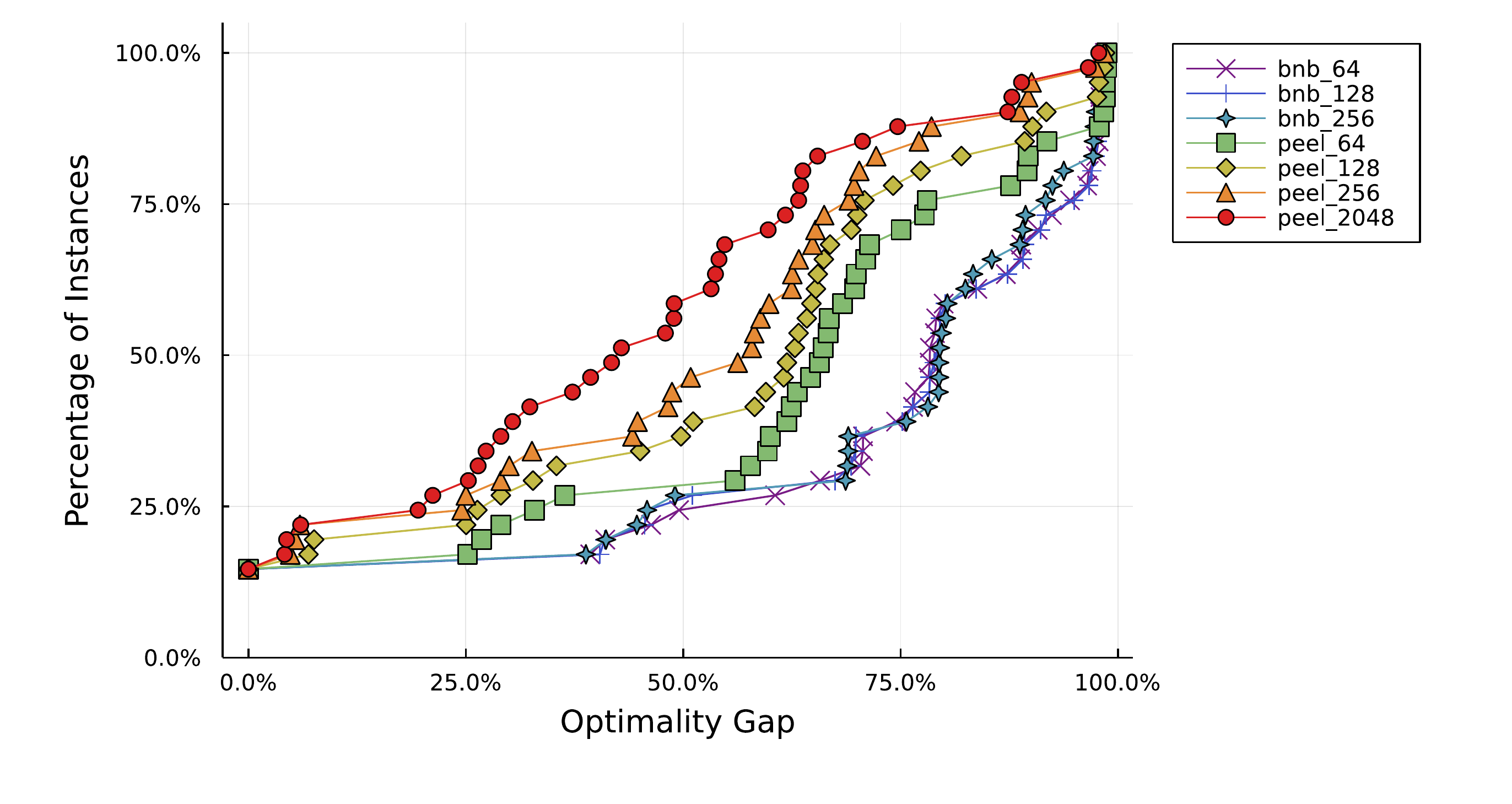}
    \caption{Performance Profiles: the optimality gap $= \frac{upper\_bound - lower\_bound}{upper\_bound}$}
    \label{fig:perfprof2}
\end{figure}

The smallest DD width tested for both methods was $64$, and the largest DD width tested was $256$. Figure \ref{fig:summary} has summary statistics for those widths as the percentage improvement demonstrated by PnB. A positive percentage always indicates that PnB performed better than BnB in that category, while a negative percentage indicates that BnB performed better. Figure \ref{fig:perfprof2} shows performance profiles for all of the experiments. Figure \ref{fig:Data4} contains summary statistics comparing PnB at width $256$ to PnB at width $2048$, where a positive percentage always indicates that the width of $2048$ performed better. 

\begin{figure}[!ht]
    \centering
    \begin{tabular}{c|c c c c | c c c c|}
        & \multicolumn{4}{c|}{Width: 64} & \multicolumn{4}{c|}{Width: 256} \\
        & RB & BS & OG & QL & RB & BS & OG & QL\\ \hline
        Average \% Improvement & $114\%$ & $0.5\%$ & $22.8\%$ & $1,647\%$ & $545\%$ & $3.3\%$ & $181\%$ & $308\%$\\
        Median \% Improvement & $26\%$ & $0.05\%$ & $17.4\%$ & $734\%$ & $80\%$ & $1.7\%$ & $35\%$ & $141\%$
    \end{tabular}
    \caption{Summary Statistics: percentage improvement of peel-and-bound over branch-and-bound. RB = Relaxed Bound, BS = Best Solution, OG = Optimality Gap, QL = Queue Length. Figures \ref{fig:Data1} and \ref{fig:Data2} in Appendix \ref{sec:data} show the comprehensive results.}
    \label{fig:summary}
\end{figure}
\begin{figure}[!ht]
    \centering
    \begin{tabular}{c|c c c |}
        & \multicolumn{3}{c|}{PnB: 2048 v PnB: 256}\\
        & Relaxed Bound & Best Solution & Optimality Gap\\ \hline
         Average \% Improvement & $19.5\%$ & $0.8\%$ & $18.6\%$ \\
         Median \% Improvement &  $16.3\%$ & $0.5\%$ & $13.7\%$
    \end{tabular}
    \caption{Summary Statistics: percentage improvement of peel-and-bound at width $2048$ over peel-and-bound at width $256$. Figure \ref{fig:Data3} in Appendix \ref{sec:data} shows the comprehensive results.}
    \label{fig:Data4}
\end{figure}

As shown in Figure~\ref{fig:summary}, peel-and-bound vastly outperforms branch-and-bound in these experiments. The average and median improvements from using peel-and-bound at both widths are significant in terms of the relaxed bound, the remaining optimality gap, and the number of nodes that still need to be processed. The best solution found by the end of the runtime also tends to be slightly better with peel-and-bound, but the found solutions are often so close to the real optimal solutions that there is little room for improvement. At both widths, six of the problems were solved to optimality. BnB was faster in only one of those cases, and in that case the difference was $.04$ seconds. The median time for PnB to close in these cases was $191\%$ faster at a width of 64, and $580\%$ faster at a width of $256$. The relaxed bound produced by PnB at a width of $64$ was better for $28$ of the remaining $35$ problems, and at a width of $256$ was better for $34$ of the remaining $35$ problems. The optimality gap was similarly better for peel-and-bound on every problem except the ones where branch-and-bound found a better relaxed bound. However, of the problems where branch-and-bound had a better optimality gap, the improvement was less than $1\%$ for all but one problem. 
Figure~\ref{fig:perfprof2} reinforces that even though there are some instances where a specific branch-and-bound setting slightly outperforms a specific peel-and-bound setting, the gap in those cases is small relative to the general gap between all peel-and-bound settings and all branch-and-bound settings. 

As shown in Figure~\ref{fig:Data4}, increasing the width to $2048$ from $256$ led to an $19.5\%$ average improvement ($16.3\%$ median improvement) in the relaxed bound. 
Figure~\ref{fig:perfprof2} shows that the performance of peel-and-bound nearly uniformly increases with the maximum allowable width. Similar to the difference between branch-and-bound and peel-and-bound, some specific instances see a small out-performance of the peel-and-bound running at a smaller width, but the gap is small relative to the usual gap between the $2048$-width experiment and the rest of the experiments. 
Additionally, Figure~\ref{fig:perfprof2} shows that peel$\_2048$ solved $50\%$ of instances to within a $42\%$ optimality gap, peel$\_64$ solved $50\%$ of instances to within a $67\%$ optimality gap, and the best performing branch and bound (bnb$\_64$) solved $50\%$ of instances to within only a $79\%$ optimality gap. The overall performance of peel-and-bound improves as more problems are considered, especially as the maximum allowable width for the decision diagrams is increased.

\begin{figure}[!ht]
    \centering
        \setlength{\tabcolsep}{1pt}
        \begin{tabular}{c c}
            \includegraphics[scale=.24]{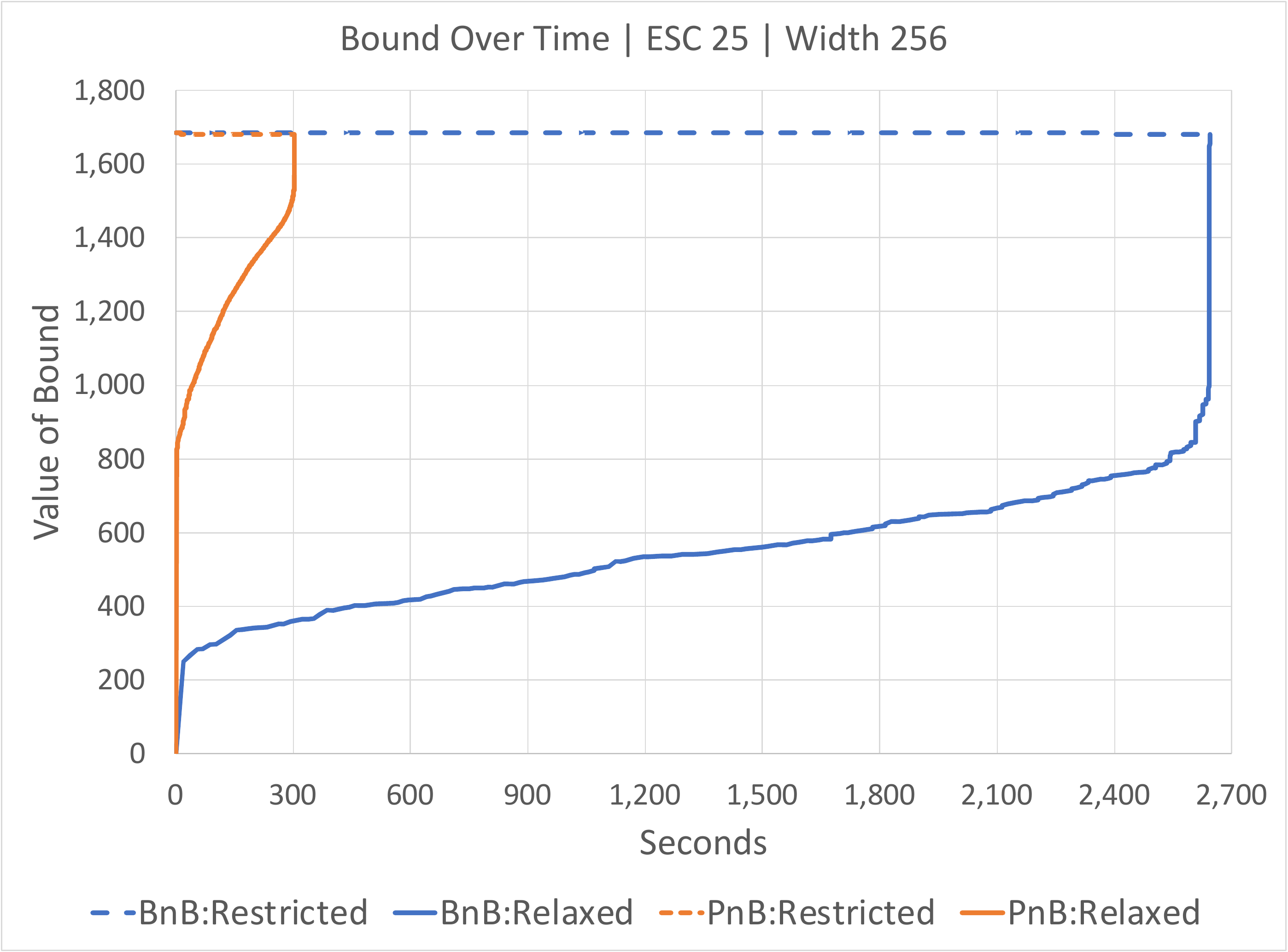} & \includegraphics[scale=.24]{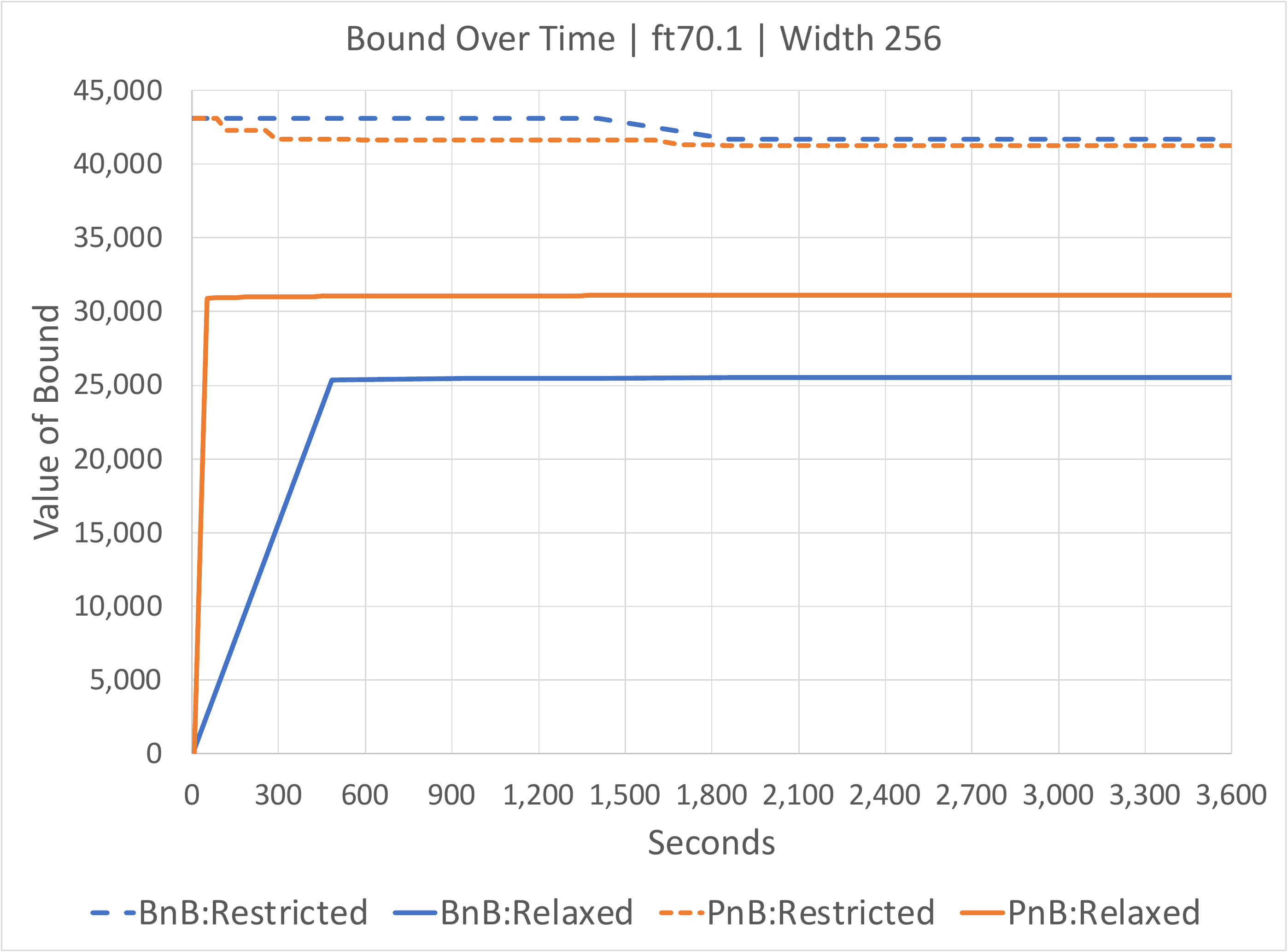}\\
           Solved SOP & Unsolved SOP
        \end{tabular}
    \caption{Dual bounds of ESC25 and ft70.1 over the runtime of the experiment.}
    \label{fig:sopExample}
\end{figure}

The selected graphs shown in Figure~\ref{fig:sopExample} are representative of the two main types of behavior observed over the problem set. On problems where the underlying relaxation method works well, the relaxed bound moves quickly towards convergence with the best found solution. On problems where the underlying relaxation does not work well, both algorithms slowly improve the relaxed bound, but PnB starts stronger as it can use exact arc values, and it maintains the advantage throughout. It is clear from the time-series data that to be competitive with cutting-edge solvers, peel-and-bound must be combined with other constraint programming propagators. However, it is also clear that peel-and-bound can have a significant edge over a propagator that generates the required decision diagrams from scratch at each iteration.

\section{Conclusion and Future Work}
\label{sec:conclusion}
This paper presented a peel-and-bound algorithm as an alternative to branch-and-bound. 
In peel-and-bound, constructed decision diagrams are repeatedly reused to avoid unnecessary computation. Additionally, peel-and-bound can be used in combination with a decision diagram structure that only admits exact arc values, to increase scalability and further reduce the amount of work needed at each iteration of the algorithm. We identified several heuristic decisions that can be used to adjust peel-and-bound, and provided insight into how the algorithm can be applied to other problems.

We compared the performance of a peel-and-bound scheme to a branch-and-bound scheme using the same DD based propagator. We tested both algorithms on the 41 instances of the SOP from TSPLIB. Results show that peel-and-bound significantly outperforms branch-and-bound on the SOP by generating substantially stronger relaxed bounds on instances that were not closed during the experiment, and reaching optimality faster when the instances were closed. This paper provides strong support for the value of re-using work in DD based solvers. Furthermore, peel-and-bound benefits from scaling the maximum allowable width. Thus, relaxed DDs that yield strong bounds at scale, but are too costly to generate iteratively, only need to be constructed once. The method detailed in this paper focused on DDs generated by separation; future research could focus on DDs generated using a merge operator.

\clearpage
\bibliography{peel-and-bound}

\appendix
\clearpage
\section{Experimental Data}
\label{sec:data}

\begin{figure}[!h]
    \centering
    \resizebox{\textwidth}{!}{%
    \begin{tabular}{|cc|ccccc|ccccc|ccccc|}
        \multicolumn{2}{c|}{Problem Info} & \multicolumn{5}{c|}{BnB: width 64} & \multicolumn{5}{c|}{PnB: width 64} & \multicolumn{5}{c}{Percent Improvements} \\
        Name & n & \multicolumn{1}{c}{RB} & \multicolumn{1}{c}{BS} & T & OG & QL & \multicolumn{1}{c}{RB} & \multicolumn{1}{c}{BS} & T & OG & QL & \multicolumn{1}{c}{RB} & BS & T & OG & QL \\ \hline
        ESC07 & 9 & 2,125 & 2,125 & 0.03 & 0\% & - & 2,125 & 2,125 & 0.07 & 0\% & - &  &  & -57\% &  &  \\
        ESC11 & 13 & 2,075 & 2,075 & 0.65 & 0\% & - & 2,075 & 2,075 & 0.42 & 0\% & - &  &  & 55\% &  &  \\
        ESC12 & 14 & 1,675 & 1,675 & 1.99 & 0\% & - & 1,675 & 1,675 & 0.64 & 0\% & - &  &  & 211\% &  &  \\
        ESC25 & 27 & 1,681 & 1,681 & 956 & 0\% & - & 1,681 & 1,681 & 353 & 0\% & - &  &  & 171\% &  &  \\
        ESC47 & 49 & 334 & 1,542 &  & 78\% & 8,842 & 368 & 1,676 &  & 78\% & 1,295 & 10.2\% & -8.0\% &  & 0.4\% & 583\% \\
        ESC63 & 65 & 8 & 62 &  & 87\% & 2,756 & 44 & 62 &  & 29\% & 15 & 450.0\% & 0.0\% &  & 200.0\% & 18273\% \\
        ESC78 & 80 & 2,230 & 19,800 &  & 89\% & 1,040 & 5,000 & 20,045 &  & 75\% & 316 & 124.2\% & -1.2\% &  & 18.2\% & 229\% \\
        br17.10 & 18 & 55 & 55 & 260 & 0\% & - & 55 & 55 & 5 & 0\% & - &  &  & 4652\% &  &  \\
        br17.12 & 18 & 55 & 55 & 138 & 0\% & - & 55 & 55 & 21 & 0\% & - &  &  & 546\% &  &  \\
        ft53.1 & 54 & 1,785 & 8,478 &  & 79\% & 8,841 & 3,324 & 8,244 &  & 60\% & 917 & 86.2\% & 2.8\% &  & 32.3\% & 864\% \\
        ft53.2 & 54 & 1,946 & 8,927 &  & 78\% & 7,356 & 3,450 & 8,633 &  & 60\% & 938 & 77.3\% & 3.4\% &  & 30.3\% & 684\% \\
        ft53.3 & 54 & 2,546 & 12,179 &  & 79\% & 5,594 & 4,234 & 12,327 &  & 66\% & 1,147 & 66.3\% & -1.2\% &  & 20.5\% & 388\% \\
        ft53.4 & 54 & 3,780 & 14,811 &  & 74\% & 11,907 & 6,500 & 14,753 &  & 56\% & 2,372 & 72.0\% & 0.4\% &  & 33.1\% & 402\% \\
        ft70.1 & 71 & 25,444 & 41,926 &  & 39\% & 4,781 & 31,123 & 41,607 &  & 25\% & 412 & 22.3\% & 0.8\% &  & 56.0\% & 1060\% \\
        ft70.2 & 71 & 25,239 & 42,805 &  & 41\% & 3,998 & 31,195 & 42,623 &  & 27\% & 427 & 23.6\% & 0.4\% &  & 53.1\% & 836\% \\
        ft70.3 & 71 & 25,810 & 48,073 &  & 46\% & 4,036 & 31,872 & 47,491 &  & 33\% & 475 & 23.5\% & 1.2\% &  & 40.8\% & 750\% \\
        ft70.4 & 71 & 28,593 & 56,644 &  & 50\% & 8,642 & 35,974 & 56,552 &  & 36\% & 1,087 & 25.8\% & 0.2\% &  & 36.1\% & 695\% \\
        kro124p.1 & 101 & 10,773 & 46,158 &  & 77\% & 2,173 & 17,579 & 46,158 &  & 62\% & 105 & 63.2\% & 0.0\% &  & 23.8\% & 1970\% \\
        kro124p.2 & 101 & 11,061 & 46,930 &  & 76\% & 1,898 & 17,633 & 46,930 &  & 62\% & 109 & 59.4\% & 0.0\% &  & 22.4\% & 1641\% \\
        kro124p.3 & 101 & 12,110 & 55,991 &  & 78\% & 1,055 & 18,586 & 55,991 &  & 67\% & 117 & 53.5\% & 0.0\% &  & 17.3\% & 802\% \\
        kro124p.4 & 101 & 13,838 & 85,533 &  & 84\% & 2,990 & 24,388 & 85,316 &  & 71\% & 244 & 76.2\% & 0.3\% &  & 17.4\% & 1125\% \\
        p43.1 & 44 & 630 & 29,450 &  & 98\% & 12,945 & 380 & 29,380 &  & 99\% & 1,022 & -39.7\% & 0.2\% &  & -0.9\% & 1167\% \\
        p43.2 & 44 & 440 & 29,000 &  & 98\% & 8,519 & 420 & 29,080 &  & 99\% & 1,125 & -4.5\% & -0.3\% &  & -0.1\% & 657\% \\
        p43.3 & 44 & 595 & 29,530 &  & 98\% & 12,802 & 490 & 29,530 &  & 98\% & 1,122 & -17.6\% & 0.0\% &  & -0.4\% & 1041\% \\
        p43.4 & 44 & 1,370 & 83,855 &  & 98\% & 21,105 & 1,050 & 83,890 &  & 99\% & 4,694 & -23.4\% & 0.0\% &  & -0.4\% & 350\% \\
        prob.42 & 42 & 99 & 289 &  & 66\% & 16,742 & 97 & 286 &  & 66\% & 2,613 & -2.0\% & 1.0\% &  & -0.5\% & 541\% \\
        prob.100 & 100 & 170 & 1,841 &  & 91\% & 1,731 & 182 & 1,760 &  & 90\% & 117 & 7.1\% & 4.6\% &  & 1.2\% & 1379\% \\
        rbg048a & 50 & 76 & 379 &  & 80\% & 12,938 & 47 & 380 &  & 88\% & 1,551 & -38.2\% & -0.3\% &  & -8.8\% & 734\% \\
        rbg050c & 52 & 63 & 566 &  & 89\% & 11,480 & 154 & 512 &  & 70\% & 1,481 & 144.4\% & 10.5\% &  & 27.1\% & 675\% \\
        rbg109a & 111 & 91 & 1,196 &  & 92\% & 2,773 & 379 & 1,196 &  & 68\% & 612 & 316.5\% & 0.0\% &  & 35.3\% & 353\% \\
        rbg150a & 152 & 63 & 1,874 &  & 97\% & 241 & 565 & 1,865 &  & 70\% & 222 & 796.8\% & 0.5\% &  & 38.6\% & 9\% \\
        rbg174a & 176 & 119 & 2,157 &  & 94\% & 809 & 626 & 2,156 &  & 71\% & 117 & 426.1\% & 0.0\% &  & 33.1\% & 591\% \\
        rbg253a & 255 & 113 & 3,181 &  & 96\% & 403 & 708 & 3,180 &  & 78\% & 39 & 526.5\% & 0.0\% &  & 24.1\% & 933\% \\
        rbg323a & 325 & 89 & 3,519 &  & 97\% & 437 & 289 & 3,529 &  & 92\% & 17 & 224.7\% & -0.3\% &  & 6.2\% & 2471\% \\
        rbg341a & 343 & 68 & 3,038 &  & 98\% & 366 & 321 & 3,064 &  & 90\% & 8 & 372.1\% & -0.8\% &  & 9.2\% & 4475\% \\
        rbg358a & 360 & 69 & 3,359 &  & 98\% & 289 & 73 & 3,373 &  & 98\% & 6 & 5.8\% & -0.4\% &  & 0.1\% & 4717\% \\
        rbg378a & 380 & 52 & 3,429 &  & 98\% & 266 & 50 & 3,429 &  & 99\% & 5 & -3.8\% & 0.0\% &  & -0.1\% & 5220\% \\
        ry48p.1 & 49 & 5,201 & 17,555 &  & 70\% & 10,480 & 6,171 & 17,454 &  & 65\% & 1,395 & 18.7\% & 0.6\% &  & 8.9\% & 651\% \\
        ry48p.2 & 49 & 5,291 & 18,046 &  & 71\% & 9,286 & 6,577 & 17,840 &  & 63\% & 1,445 & 24.3\% & 1.2\% &  & 12.0\% & 543\% \\
        ry48p.3 & 49 & 6,207 & 21,161 &  & 71\% & 9,039 & 6,985 & 20,962 &  & 67\% & 1,707 & 12.5\% & 0.9\% &  & 6.0\% & 430\% \\
        ry48p.4 & 49 & 13,610 & 34,517 &  & 61\% & 15,819 & 14,293 & 33,804 &  & 58\% & 3,217 & 5.0\% & 2.1\% &  & 4.9\% & 392\%
    \end{tabular}
    }
     \caption{Comparison Data for width 64 experiments: RB = Relaxed Bound, BS = Best Solution, T = Time in Seconds, OG = Optimality Gap, QL = Queue Length. Full time series data is available in the GitHub repository.}
     \label{fig:Data1}
\end{figure}

\begin{figure}[!ht]
    \centering
    \resizebox{\textwidth}{!}{%
    \begin{tabular}{|cc|ccccc|ccccc|ccccc}
        \multicolumn{2}{c|}{Problem   Info} & \multicolumn{5}{c|}{BnB: width 256} & \multicolumn{5}{c|}{PnB: width 256} & \multicolumn{5}{c}{Percent Improvements} \\
        Name & n & \multicolumn{1}{c}{RB} & \multicolumn{1}{c}{BS} & T & OG & QL & \multicolumn{1}{c}{RB} & \multicolumn{1}{c}{BS} & T & OG & QL & \multicolumn{1}{c}{RB} & BS & T & OG & QL \\ \hline
        ESC07 & 9 & 2,125 & 2,125 & 0.04 & 0\% & - & 2,125 & 2,125 & 0.04 & 0\% & - &  &  & 0\% &  &  \\
        ESC11 & 13 & 2,075 & 2,075 & 0.48 & 0\% & - & 2,075 & 2,075 & 0.41 & 0\% & - &  &  & 17\% &  &  \\
        ESC12 & 14 & 1,675 & 1,675 & 1.66 & 0\% & - & 1,675 & 1,675 & 0.34 & 0\% & - &  &  & 388\% &  &  \\
        ESC25 & 27 & 1,681 & 1,681 & 2,643 & 0\% & - & 1,681 & 1,681 & 303 & 0\% & - &  &  & 771\% &  &  \\
        ESC47 & 49 & 312 & 1,590 &  & 80\% & 720 & 658 & 1,339 &  & 51\% & 740 & 110.9\% & 18.7\% &  & 58.0\% & -3\% \\
        ESC63 & 65 & 9 & 62 &  & 85\% & 53 & 44 & 62 &  & 29\% & 3 & 388.9\% & 0.0\% &  & 194.4\% & 1667\% \\
        ESC78 & 80 & 2,230 & 20,345 &  & 89\% & 59 & 5,600 & 20,135 &  & 72\% & 109 & 151.1\% & 1.0\% &  & 23.3\% & -46\% \\
        br17.10 & 18 & 55 & 55 & 275 & 0\% & - & 55 & 55 & 3 & 0\% & - &  &  & 9468\% &  &  \\
        br17.12 & 18 & 55 & 55 & 105 & 0\% & - & 55 & 55 & 5 & 0\% & - &  &  & 2146\% &  &  \\
        ft53.1 & 54 & 1,708 & 8,424 &  & 80\% & 760 & 4,603 & 8,244 &  & 44\% & 271 & 169.5\% & 2.2\% &  & 80.5\% & 180\% \\
        ft53.2 & 54 & 1,856 & 9,059 &  & 80\% & 632 & 3,555 & 8,648 &  & 59\% & 272 & 91.5\% & 4.8\% &  & 35.0\% & 132\% \\
        ft53.3 & 54 & 2,493 & 12,598 &  & 80\% & 477 & 4,852 & 11,095 &  & 56\% & 390 & 94.6\% & 13.5\% &  & 42.6\% & 22\% \\
        ft53.4 & 54 & 3,619 & 14,867 &  & 76\% & 1,240 & 7,560 & 14,611 &  & 48\% & 797 & 108.9\% & 1.8\% &  & 56.8\% & 56\% \\
        ft70.1 & 71 & 25,507 & 41,686 &  & 39\% & 373 & 31,122 & 41,235 &  & 25\% & 108 & 22.0\% & 1.1\% &  & 58.3\% & 245\% \\
        ft70.2 & 71 & 25,261 & 42,901 &  & 41\% & 297 & 31,630 & 42,182 &  & 25\% & 123 & 25.2\% & 1.7\% &  & 64.4\% & 141\% \\
        ft70.3 & 71 & 25,891 & 47,806 &  & 46\% & 377 & 32,539 & 46,488 &  & 30\% & 151 & 25.7\% & 2.8\% &  & 52.8\% & 150\% \\
        ft70.4 & 71 & 31,186 & 56,366 &  & 45\% & 958 & 37,984 & 56,366 &  & 33\% & 356 & 21.8\% & 0.0\% &  & 37.0\% & 169\% \\
        kro124p.1 & 101 & 10,683 & 48,866 &  & 78\% & 152 & 19,224 & 45,643 &  & 58\% & 43 & 79.9\% & 7.1\% &  & 35.0\% & 253\% \\
        kro124p.2 & 101 & 10,706 & 52,038 &  & 79\% & 125 & 19,299 & 48,102 &  & 60\% & 43 & 80.3\% & 8.2\% &  & 32.6\% & 191\% \\
        kro124p.3 & 101 & 12,078 & 58,562 &  & 79\% & 64 & 20,145 & 57,358 &  & 65\% & 45 & 66.8\% & 2.1\% &  & 22.3\% & 42\% \\
        kro124p.4 & 101 & 14,511 & 82,672 &  & 82\% & 281 & 25,002 & 82,364 &  & 70\% & 102 & 72.3\% & 0.4\% &  & 18.4\% & 175\% \\
        p43.1 & 44 & 610 & 29,460 &  & 98\% & 1,033 & 27,255 & 28,635 &  & 5\% & 146 & 4368\% & 2.9\% &  & 1932\% & 608\% \\
        p43.2 & 44 & 460 & 29,020 &  & 98\% & 547 & 27,455 & 29,020 &  & 5\% & 391 & 5868\% & 0.0\% &  & 1725\% & 40\% \\
        p43.3 & 44 & 750 & 29,530 &  & 97\% & 1,016 & 27,780 & 29,530 &  & 6\% & 764 & 3604\% & 0.0\% &  & 1545\% & 33\% \\
        p43.4 & 44 & 1,425 & 83,880 &  & 98\% & 1,365 & 28,195 & 83,435 &  & 66\% & 1,380 & 1879\% & 0.5\% &  & 48.5\% & -1\% \\
        prob.42 & 42 & 90 & 289 &  & 69\% & 1,166 & 103 & 275 &  & 63\% & 617 & 14.4\% & 5.1\% &  & 10.1\% & 89\% \\
        prob.100 & 100 & 157 & 1,886 &  & 92\% & 113 & 178 & 1,721 &  & 90\% & 45 & 13.4\% & 9.6\% &  & 2.3\% & 151\% \\
        rbg048a & 50 & 80 & 389 &  & 79\% & 794 & 80 & 373 &  & 79\% & 534 & 0.0\% & 4.3\% &  & 1.1\% & 49\% \\
        rbg050c & 52 & 62 & 583 &  & 89\% & 810 & 175 & 503 &  & 65\% & 442 & 182.3\% & 15.9\% &  & 37.0\% & 83\% \\
        rbg109a & 111 & 89 & 1,181 &  & 92\% & 394 & 406 & 1,106 &  & 63\% & 204 & 356.2\% & 6.8\% &  & 46.1\% & 93\% \\
        rbg150a & 152 & 115 & 1,845 &  & 94\% & 406 & 571 & 1,845 &  & 69\% & 100 & 396.5\% & 0.0\% &  & 35.8\% & 306\% \\
        rbg174a & 176 & 362 & 2,172 &  & 83\% & 337 & 646 & 2,171 &  & 70\% & 57 & 78.5\% & 0.0\% &  & 18.6\% & 491\% \\
        rbg253a & 255 & 359 & 3,177 &  & 89\% & 139 & 727 & 3,176 &  & 77\% & 22 & 102.5\% & 0.0\% &  & 15.0\% & 532\% \\
        rbg323a & 325 & 99 & 3,476 &  & 97\% & 114 & 346 & 3,480 &  & 90\% & 14 & 249.5\% & -0.1\% &  & 7.9\% & 714\% \\
        rbg341a & 343 & 84 & 3,016 &  & 97\% & 120 & 340 & 3,016 &  & 89\% & 7 & 304.8\% & 0.0\% &  & 9.6\% & 1614\% \\
        rbg358a & 360 & 88 & 3,280 &  & 97\% & 92 & 88 & 3,382 &  & 97\% & 5 & 0.0\% & -3.0\% &  & -0.1\% & 1740\% \\
        rbg378a & 380 & 44 & 3,385 &  & 99\% & 35 & 53 & 3,385 &  & 98\% & 6 & 20.5\% & 0.0\% &  & 0.3\% & 483\% \\
        ry48p.1 & 49 & 5,470 & 17,464 &  & 69\% & 897 & 9,432 & 17,071 &  & 45\% & 377 & 72.4\% & 2.3\% &  & 53.5\% & 138\% \\
        ry48p.2 & 49 & 5,606 & 18,060 &  & 69\% & 834 & 6,615 & 17,627 &  & 62\% & 383 & 18.0\% & 2.5\% &  & 10.4\% & 118\% \\
        ry48p.3 & 49 & 6,558 & 21,142 &  & 69\% & 859 & 8,723 & 20,850 &  & 58\% & 513 & 33.0\% & 1.4\% &  & 18.6\% & 67\% \\
        ry48p.4 & 49 & 17,359 & 34,074 &  & 49\% & 1,557 & 17,322 & 33,773 &  & 49\% & 990 & -0.2\% & 0.9\% &  & 0.7\% & 57\%
    \end{tabular}
    }
     \caption{Comparison data for width 256 experiments: RB = Relaxed Bound, BS = Best Solution, T = Time in Seconds, OG = Optimality Gap, QL = Queue Length. Full time series data is available in the GitHub repository.}
     \label{fig:Data2}
 \end{figure}

\begin{figure}[!ht]
    \centering
    \resizebox{\textwidth}{!}{%
    \begin{tabular}{|cc|ccc|ccc|ccc|}
        \multicolumn{2}{|c|}{Problem   Info} & \multicolumn{3}{c|}{PnB: width 256} & \multicolumn{3}{c|}{PnB: width 2048} & \multicolumn{3}{c|}{Percent Improvements} \\
        Name & n & \multicolumn{1}{c}{RB} & \multicolumn{1}{c}{BS} & OG & \multicolumn{1}{c}{RB} & \multicolumn{1}{c}{BS} & OG & \multicolumn{1}{c}{RB} & BS & OG \\ \hline
        ESC47 & 49 & 658 & 1,339 & 51\% & 882 & 1,304 & 32\% & 34.0\% & 2.7\% & 57.2\% \\
        ESC63 & 65 & 44 & 62 & 29\% & 44 & 62 & 29\% & 0.0\% & 0.0\% & 0\% \\
        ESC78 & 80 & 5,600 & 20,135 & 72\% & 6,025 & 20,505 & 71\% & 7.6\% & -1.8\% & 2.2\% \\
        ft53.1 & 54 & 4,603 & 8,244 & 44\% & 5,167 & 8,237 & 37\% & 12.3\% & 0.1\% & 18.5\% \\
        ft53.2 & 54 & 3,555 & 8,648 & 59\% & 4,910 & 8,598 & 43\% & 38.1\% & 0.6\% & 37.3\% \\
        ft53.3 & 54 & 4,852 & 11,095 & 56\% & 7,722 & 11,092 & 30\% & 59.2\% & 0.0\% & 85.2\% \\
        ft53.4 & 54 & 7,560 & 14,611 & 48\% & 7,466 & 14,618 & 49\% & -1.2\% & 0.0\% & -1.4\% \\
        ft70.1 & 71 & 31,122 & 41,235 & 25\% & 33,382 & 41,476 & 20\% & 7.3\% & -0.6\% & 25.7\% \\
        ft70.2 & 71 & 31,630 & 42,182 & 25\% & 32,964 & 41,833 & 21\% & 4.2\% & 0.8\% & 18.0\% \\
        ft70.3 & 71 & 32,539 & 46,488 & 30\% & 34,366 & 46,001 & 25\% & 5.6\% & 1.1\% & 18.6\% \\
        ft70.4 & 71 & 37,984 & 56,366 & 33\% & 40,919 & 56,310 & 27\% & 7.7\% & 0.1\% & 19.3\% \\
        kro124p.1 & 101 & 19,224 & 45,643 & 58\% & 21,954 & 47,425 & 54\% & 14.2\% & -3.8\% & 7.8\% \\
        kro124p.2 & 101 & 19,299 & 48,102 & 60\% & 22,746 & 49,571 & 54\% & 17.9\% & -3.0\% & 10.7\% \\
        kro124p.3 & 101 & 20,145 & 57,358 & 65\% & 25,566 & 54,633 & 53\% & 26.9\% & 5.0\% & 21.9\% \\
        kro124p.4 & 101 & 25,002 & 82,364 & 70\% & 29,377 & 81,050 & 64\% & 17.5\% & 1.6\% & 9.2\% \\
        p43.1 & 44 & 27,255 & 28,635 & 5\% & 27,755 & 28,960 & 4\% & 1.8\% & -1.1\% & 16\% \\
        p43.2 & 44 & 27,455 & 29,020 & 5\% & 27,725 & 29,000 & 4\% & 1.0\% & 0.1\% & 23\% \\
        p43.3 & 44 & 27,780 & 29,530 & 6\% & 27,755 & 29,530 & 6\% & -0.1\% & 0.0\% & -1\% \\
        p43.4 & 44 & 28,195 & 83,435 & 66\% & 28,680 & 83,020 & 65\% & 1.7\% & 0.5\% & 1.2\% \\
        prob.42 & 42 & 103 & 275 & 63\% & 152 & 261 & 42\% & 47.6\% & 5.4\% & 49.8\% \\
        prob.100 & 100 & 178 & 1,721 & 90\% & 220 & 1,735 & 87\% & 23.6\% & -0.8\% & 2.7\% \\
        rbg048a & 50 & 80 & 373 & 79\% & 93 & 367 & 75\% & 16.3\% & 1.6\% & 5.2\% \\
        rbg050c & 52 & 175 & 503 & 65\% & 184 & 501 & 63\% & 5.1\% & 0.4\% & 3.1\% \\
        rbg109a & 111 & 406 & 1,106 & 63\% & 453 & 1,126 & 60\% & 11.6\% & -1.8\% & 5.9\% \\
        rbg150a & 152 & 571 & 1,845 & 69\% & 672 & 1,841 & 63\% & 17.7\% & 0.2\% & 8.7\% \\
        rbg174a & 176 & 646 & 2,171 & 70\% & 1,104 & 2,121 & 48\% & 70.9\% & 2.4\% & 46.5\% \\
        rbg253a & 255 & 727 & 3,176 & 77\% & 1,186 & 3,101 & 62\% & 63.1\% & 2.4\% & 24.9\% \\
        rbg323a & 325 & 346 & 3,480 & 90\% & 421 & 3,449 & 88\% & 21.7\% & 0.9\% & 2.6\% \\
        rbg341a & 343 & 340 & 3,016 & 89\% & 329 & 2,965 & 89\% & -3.2\% & 1.7\% & -0.2\% \\
        rbg358a & 360 & 88 & 3,382 & 97\% & 107 & 3,131 & 97\% & 21.6\% & 8.0\% & 0.8\% \\
        rbg378a & 380 & 53 & 3,385 & 98\% & 74 & 3,338 & 98\% & 39.6\% & 1.4\% & 0.7\% \\
        ry48p.1 & 49 & 9,432 & 17,071 & 45\% & 10,386 & 17,124 & 39\% & 10.1\% & -0.3\% & 13.7\% \\
        ry48p.2 & 49 & 6,615 & 17,627 & 62\% & 7,896 & 17,461 & 55\% & 19.4\% & 1.0\% & 14.0\% \\
        ry48p.3 & 49 & 8,723 & 20,850 & 58\% & 10,558 & 20,686 & 49\% & 21.0\% & 0.8\% & 18.8\% \\
        ry48p.4 & 49 & 17,322 & 33,773 & 49\% & 24,248 & 32,953 & 26\% & 40.0\% & 2.5\% & 84.4\%
    \end{tabular}
    }
    \caption{Comparison of PnB at 2048 over PnB at 256: RB = Relaxed Bound, BS = Best Solution, OG = Optimality Gap.}
    \label{fig:Data3}
 \end{figure}
\end{document}